\documentclass[a4paper,11pt,leqno]{article}%
\usepackage{amsmath, amssymb, amsthm}
\usepackage{graphicx}
\usepackage[latin1]{inputenc}
\usepackage[all]{xy}
\usepackage{amsmath}
\usepackage{amsfonts}
\usepackage{amssymb}
\usepackage[active]{srcltx}
\usepackage[normalem]{ulem}%
\providecommand{\U}[1]{\protect\rule{.1in}{.1in}}

\def\MM{{\mathcal M}}

\def\SSS{{\mathfrak S}}
\def\Z{{\mathbb Z}}

\def\F{{\mathbb F}}

\def\bb{{\mathbf b}}

\providecommand{\U}[1]{\protect\rule{.1in}{.1in}}
\nonstopmode

\def\Changed/{\ifvmode\else\vadjust{\vbox to 0pt{\vskip -\baselineskip\hbox to 0pt{\hss\vrule height 0pt depth 1.2\baselineskip\hskip 1em}\vss}}\fi}
\def\Math#1{\def\MathString{#1}\futurelet\MathDelim\MathChoose}
\def\MathChoose{\ifmmode\let\MathDo\MathString              \else\let\MathDo\MathSkip\fi              \MathDo}
\def\MathSkip{\ifx\MathDelim/\def\MathDo{$\MathString$\EatOne}              \else\def\MathDo{$\MathString$}\fi              \MathDo}
\def\Text#1{\def\TextString{#1}\futurelet\TextDelim\TextSkip}
\def\TextSkip{\ifx\TextDelim/\def\TextDo{\TextString\EatOne}              \else\let\TextDo\TextString\fi              \TextDo}
\def\EatOne#1{}
\def\SkipToEndScan#1\EndScan{}
\def\Scan#1#2#3{\ifx#1#2#3\expandafter\SkipToEndScan\fi\Scan#1}
\def\Upper#1{\Scan#1aAbBcCdDeEfFgGhHiIjJkKlLmMnNoOpPqQrRsStTuUvVwWxXyYzZ#1#1\EndScan}
\def\Phrase#1 #2/#3/#4=#5 #6/#7/#8.{\expandafter\edef\csname#2#3\endcsname{\noexpand\Text{#6#7}}
\expandafter\edef\csname\Upper#2#3\endcsname{\noexpand\Text{\Upper#6#7}}
\expandafter\edef\csname#1#2#3\endcsname{\noexpand\Text{#5 #6#7}}
\expandafter\edef\csname\Upper#1#2#3\endcsname{\noexpand\Text{\Upper#5 #6#7}}
\expandafter\edef\csname#2#4\endcsname{\noexpand\Text{#6#8}}
\expandafter\edef\csname\Upper#2#4\endcsname{\noexpand\Text{\Upper#6#8}}
}

\theoremstyle{plain}
\newtheorem{Theorem}{Theorem}[section]
\newtheorem{theorem}[Theorem]{Theorem}

\newtheorem{corollary}[Theorem]{Corollary}

\newtheorem{lemma}[Theorem]{Lemma}

\newtheorem{proposition}[Theorem]{Proposition}
\theoremstyle{definition}

\begin{document}

\title{\textbf{{Finite index subgroups of mapping class groups}}}
\author{\textsc{A.\,J.\ Berrick, V.\ Gebhardt and L.\ Paris}}
\maketitle

\begin{abstract}
\noindent Let $g\geq3$ and $n\geq0$, and let ${\mathcal{M}}_{g,n}$ be the
mapping class group of a surface of genus $g$ with $n$ boundary components. We
prove that ${\mathcal{M}}_{g,n}$ contains a unique subgroup of index
$2^{g-1}(2^{g}-1)$ up to conjugation, a unique subgroup of index
$2^{g-1}(2^{g}+1)$ up to conjugation, and the other proper subgroups of
${\mathcal{M}}_{g,n}$ are of index greater than $2^{g-1}(2^{g}+1)$. In
particular, the minimum index for a proper subgroup of ${\mathcal{M}}_{g,n}$
is $2^{g-1}(2^{g}-1)$.

\end{abstract}

\numberwithin{equation}{section}




\pagestyle{myheadings}
\markright{{Finite index subgroups of mapping class groups}  \quad\sc{\today}}
\setcounter{section}{-1}%


\section{Introduction and statement of results\label{11107S0}}

The interaction between mapping class groups and finite groups has long been a
topic of interest. The famous Hurwitz bound of 1893 showed that a closed
Riemann surface of genus $g$ has an upper bound of $84(g-1)$ for the order of
its finite subgroups, and Kerckhoff showed that the order of finite cyclic
subgroups is bounded above by $4g+2$ \cite{Hurwi1}, \cite{Kerko1}.

The subject of finite index subgroups of mapping class groups was brought into
focus by Grossman's discovery that the mapping class group ${\mathcal{M}%
}_{g,n}=\pi_{0}(\mathrm{Homeo}(\Sigma_{g,n}))$ of an oriented surface
$\Sigma_{g,n}$ of genus $g$ and $n$ boundary components is residually finite,
and thus well-endowed with subgroups of finite index \cite{Gross1}.
($\mathrm{Homeo}(\Sigma_{g,n})$ denotes the space of those homeomorphisms of
$\Sigma_{g,n}$ that preserve the orientation and are the identity on the
boundary.) This prompts the \textquotedblleft dual\textquotedblright%
\ question: --\smallskip

\emph{What is the minimum index } $\mathrm{mi}({\mathcal{M}}_{g,n})$ \emph{of
a proper subgroup of finite index in }${\mathcal{M}}_{g,n}$\thinspace?

\smallskip Results to date have suggested that, like the maximum finite order
question, the minimum index question should have an answer that is linear in
$g$. The best previously published bound is $\mathrm{mi}({\mathcal{M}}%
_{g,n})>4g+4$ for $g\geq3$ (see \cite{Paris1}). This inequality is used by
Aramayona and Souto to prove that, if $g\geq6$ and $g^{\prime}\leq2g-1$, then
any nontrivial homomorphism ${\mathcal{M}}_{g,n}\rightarrow{\mathcal{M}%
}_{g^{\prime},n^{\prime}}$ is induced by an embedding \cite{AraSou1}. It is
also an important ingredient in the proof of Zimmermann \cite{Zimme1} that,
for $g=3$ and $4$, the minimal nontrivial quotient of ${\mathcal{M}}_{g,0}$ is
$\mathrm{Sp}_{2g}({\mathbb{F}}_{2})$.

\medskip The \textquotedblleft headline\textquotedblright\ result of this
paper is the following exact, exponential bound.

\begin{theorem}
\label{11107C0.4} For $g\geq3$ and $n\geq0$,
\[
\mathrm{mi}({\mathcal{M}}_{g,n})=\mathrm{mi}(\mathrm{Sp}_{2g}({\mathbb{Z}%
}))=\mathrm{mi}(\mathrm{Sp}_{2g}({\mathbb{F}}_{2}))=2^{g-1}(2^{g}-1)\,.
\]

\end{theorem}

This exponential bound is all the more surprising since in similar questions
we get linear (expected) bounds. For instance, Bridson \cite{Brids1,Brids2}
has proved that a mapping class group of a surface of genus $g$ cannot act by
semisimple isometries, without a global fixed point, on a CAT(0) space of
dimension less than $g$. The exact minimal dimension for such an action is
unknown. On the other hand, it has been also shown by Bridson (see
\cite{Brids1}) that ${\mathcal{M}}_{g,n}$ has only finitely many irreducible
linear representations over any algebraically closed field, up to dimension
$(g+1)$. Later, Funar \cite{Funar1} showed that there is no linear
representation with infinite image up to dimension about $\sqrt{g+1}$.
However, there is an obvious linear representation of rank $2g$ which comes
from the action of ${\mathcal{M}}_{g,n}$ on the homology of $\Sigma_{g,0}$
(the map $\theta_{g,n}$ defined below). It is expected that this
representation is minimal in some sense (see \cite{Farb1}).

The nontrivial quotient of ${\mathcal{M}}_{g,n}$ of minimal order is unknown,
but obviously its order must be at least $\mathrm{mi}({\mathcal{M}}_{g,n})$.
This quotient is known to be $\mathrm{Sp}_{2g}({\mathbb{F}}_{2})$ if
$(g,n)=(3,0)$ or $(4,0)$ (see \cite{Zimme1}). A consequence of the above
theorem is that ${\mathcal{M}}_{g^{\prime},n^{\prime}}$ cannot be a quotient
of ${\mathcal{M}}_{g,n}$ if $3\leq g^{\prime}<g$.

\medskip Our proof of this result is constructive, in ways that we now
describe. From the surface $\Sigma_{g,n}$ we obtain a closed oriented surface
$\widehat{\Sigma}_{g}$ of genus $g$ by gluing a disk along each boundary
component. The embedding $\Sigma_{g,n}\hookrightarrow\widehat{\Sigma}_{g}$
induces a first epimorphism ${\mathcal{M}}_{g,n}\twoheadrightarrow
{\mathcal{M}}_{g,0}$. The action of $\mathrm{Homeo}(\widehat{\Sigma}_{g})$ on
$H_{1}(\widehat{\Sigma}_{g})={\mathbb{Z}}^{2g}$ induces a second epimorphism
${\mathcal{M}}_{g,0}\twoheadrightarrow\mathrm{Sp}_{2g}({\mathbb{Z}})$ onto the
integral symplectic group, and, passing mod $2$, we obtain a third epimorphism
$\mathrm{Sp}_{2g}({\mathbb{Z}})\twoheadrightarrow\mathrm{Sp}_{2g}({\mathbb{F}%
}_{2})$, where ${\mathbb{F}}_{2}={\mathbb{Z}}/2{\mathbb{Z}}$. From now on we
denote by $\theta_{g,n}:{\mathcal{M}}_{g,n}\rightarrow\mathrm{Sp}%
_{2g}({\mathbb{F}}_{2})$ the composition of these three epimorphisms.

\medskip The orthogonal groups $O_{2g}^{+}({\mathbb{F}}_{2})$ and $O_{2g}%
^{-}({\mathbb{F}}_{2})$ are subgroups of $\mathrm{Sp}_{2g}({\mathbb{F}}_{2})$.
The cardinalities of $\mathrm{Sp}_{2g}({\mathbb{F}}_{2})$, $O_{2g}%
^{+}({\mathbb{F}}_{2})$ and $O_{2g}^{-}({\mathbb{F}}_{2})$ can be found for
instance in \cite{Taylo1}, and from this data it is easily shown that, for
$g\geq2$, the indices of $O_{2g}^{+}({\mathbb{F}}_{2})$ and $O_{2g}%
^{-}({\mathbb{F}}_{2})$ in $\mathrm{Sp}_{2g}({\mathbb{F}}_{2})$ are $N_{g}%
^{+}=2^{g-1}(2^{g}+1)$ and $N_{g}^{-}=2^{g-1}(2^{g}-1)$, respectively. The
following, more or less known to experts but seemingly unpublished, is the
starting-point for our main result (Theorem \ref{11107T0.2}).

\begin{theorem}
\label{11107T0.1}Let $g\geq3$.

\begin{enumerate}
\item $O_{2g}^{-}({\mathbb{F}}_{2})$ is the unique subgroup of $\mathrm{Sp}%
_{2g}({\mathbb{F}}_{2})$ of index $N_{g}^{-}$, up to conjugation.

\item $O_{2g}^{+}({\mathbb{F}}_{2})$ is the unique subgroup of $\mathrm{Sp}%
_{2g}({\mathbb{F}}_{2})$ of index $N_{g}^{+}$, up to conjugation.

\item All the other proper subgroups of $\mathrm{Sp}_{2g}({\mathbb{F}}_{2})$
are of index at least $2N_{g}^{-}$.
\end{enumerate}
\end{theorem}

\noindent We set ${\mathcal{O}}_{g,n}^{+}=\theta_{g,n}^{-1}(O_{2g}%
^{+}({\mathbb{F}}_{2}))$ and ${\mathcal{O}}_{g,n}^{-}=\theta_{g,n}^{-1}%
(O_{2g}^{-}({\mathbb{F}}_{2}))$. Thus, by the above, ${\mathcal{O}}_{g,n}^{-}$
is an index $N_{g}^{-}$ subgroup of ${\mathcal{M}}_{g,n}$, and ${\mathcal{O}%
}_{g,n}^{+}$ is an index $N_{g}^{+}$ subgroup of ${\mathcal{M}}_{g,n}$. Here
is our main result.

\begin{theorem}
\label{11107T0.2}Let $g\geq3$ and $n\geq0$.

\begin{enumerate}
\item $\mathcal{O}_{g,n}^{-}$ is the unique subgroup of ${\mathcal{M}}_{g,n}$
of index $N_{g}^{-}$, up to conjugation.

\item $\mathcal{O}_{g,n}^{+}$ is the unique subgroup of ${\mathcal{M}}_{g,n}$
of index $N_{g}^{+}$, up to conjugation.

\item If $g=3$ then all the other proper subgroups of ${\mathcal{M}}_{3,n}$
are of index strictly greater than $N_{g}^{+}=36$.

\item If $g\geq4$ then all the other proper subgroups of ${\mathcal{M}}_{g,n}$
are of index at least $5N_{g-1}^{-}>N_{g}^{+}$.
\end{enumerate}
\end{theorem}

Since $m$ is an upper bound for the minimum index of a group $G$ if and only
if there is a nontrivial homomorphism from $G$ to the symmetric group
$\mathfrak{S}_{m}$ on $m$ letters, one would like to understand the
permutation representations associated to Theorem \ref{11107T0.2}.

We denote by $\phi_{g,n}^{+}:{\mathcal{M}}_{g,n}\rightarrow{\mathfrak{S}%
}_{N_{g}^{+}}$ (resp.~$\phi_{g,n}^{-}:{\mathcal{M}}_{g,n}\rightarrow
{\mathfrak{S}}_{N_{g}^{-}}$) the permutation representation induced by the
action of ${\mathcal{M}}_{g,n}$ on the right cosets of ${\mathcal{O}}%
_{g,n}^{+}$ (resp.~${\mathcal{O}}_{g,n}^{-}$). Corresponding to the numerical
relations%
\[
N_{g}^{+}=3N_{g-1}^{+}+N_{g-1}^{-}\,,\quad N_{g}^{-}=3N_{g-1}^{-}+N_{g-1}%
^{+}\,,
\]
we prove the following.

\begin{theorem}
\label{11107T0.3}

\begin{enumerate}
\item Let $g\geq3$ and $n\geq1$. Then $\phi_{g,n}^{-}:\mathcal{M}%
_{g,n}\rightarrow\mathfrak{S}_{N_{g}^{-}}$ is, up to equivalence, the unique
extension of the representation $(\phi_{g-1,n}^{-})^{3}\oplus\phi_{g-1,n}^{+}$
from $\mathcal{M}_{g-1,n}$ to $\mathcal{M}_{g,n}$, and $\phi_{g,n}%
^{+}:\mathcal{M}_{g,n}\rightarrow\mathfrak{S}_{N_{g}^{+}}$ is, up to
equivalence, the unique extension of the representation $\phi_{g-1,n}%
^{-}\oplus(\phi_{g-1,n}^{+})^{3}$ from $\mathcal{M}_{g-1,n}$ to $\mathcal{M}%
_{g,n}$.

\item Let $g\geq3$ and $n\geq0$. Let $b$ be a nonseparating simple closed
curve on $\Sigma_{g,n}$, and let $T_{b}$ be the Dehn twist around $b$. Then
the cycle structure of the image of $T_{b}$ under $\phi_{g,n}^{-}$ is
\[
(1)^{2^{2g-2}}(2)^{2^{g-2}(2^{g-1}-1)}\,,
\]
and the cycle structure of the image of $T_{b}$ under $\phi_{g,n}^{+}$ is
\[
(1)^{2^{2g-2}}(2)^{2^{g-2}(2^{g-1}+1)}\,.
\]

\end{enumerate}
\end{theorem}

\medskip\textbf{Remark.} Implicit in the statement of Theorem \ref{11107T0.3}
is the fact that, if $n\geq1$, then ${\mathcal{M}}_{g-1,n}$ naturally embeds
into ${\mathcal{M}}_{g,n}$. This embedding will be described in Section
\ref{11107S2}. However, there is no natural embedding of ${\mathcal{M}%
}_{g-1,0}$ into ${\mathcal{M}}_{g,0}$, hence Part (1) of the theorem would
make no sense for $n=0$.

\medskip We observe that the abelianization of ${\mathcal{M}}_{g,n}$ is
isomorphic to ${\mathbb{Z}}/12\,{\mathbb{Z}}$ if $(g,n)=(1,0)$, ${\mathbb{Z}%
}^{n}$ if $g=1$ and $n\geq1$, and ${\mathbb{Z}}/10\,{\mathbb{Z}}$ if $g=2$
(see \cite{Korkm1}). Hence, the minimum index of ${\mathcal{M}}_{g,n}$ is $2$
if $g=1$ or $2$. Note that ${\mathcal{M}}_{g,n}$ is perfect if $g\geq3$ (see
\cite{Powel1, Korkm1}). If $g=2$ we have $N_{g}^{-}=6$ and $N_{g}^{+}=10$, and
there are six proper subgroups of index at most $10$ in ${\mathcal{M}}_{g,n}$:
one of index $2$, one of index $5$, two of index $6$ and two of index $10$.
The description of these subgroups as well as the proof of this fact are given
in Section \ref{11107S3}.

\medskip The remainder of the paper is divided into two parts. Part I starts
with some preliminaries on permutations (Section~\ref{11107S1}) and
presentations of mapping class groups (Section~\ref{11107S2}). Then we
determine the subgroups of ${\mathcal{M}}_{2,n}$ of index at most
$10=N_{2}^{+}$ (Section~\ref{11107S3}) and the subgroups of ${\mathcal{M}%
}_{3,n}$ of index at most $36=N_{3}^{+}$ (Section~\ref{11107S4}). We prove our
theorems by induction on the genus. The starting case, $g=3$, is made in
Section 4, and the inductive step is the object of Part II. We first treat the
case of a surface with a unique boundary component (Sections~\ref{11107S5}
to~\ref{11107S8}) and we extend the result to surfaces with several boundary
components in Section~\ref{11107S9}. Theorem 0.1 is proved in
Section~~\ref{11107S7}. This can be read independently from the rest.

\bigskip\textbf{Acknowledgements.} Many thanks to Derek Holt for his helpful
advice. The authors are pleased to acknowledge the financial support of
National University of Singapore research grants R-146-000-097-112 and
R-146-000-137-112 towards the visits of LP to Singapore in 2007 and 2011 and
AJB to Dijon in 2009 and 2010, and the warm hospitality of their hosts. LP is
partially supported by the \textit{Agence Nationale de la Recherche}
(\textit{projet Théorie de Garside}, ANR-08-BLAN-0269-03). VG is partially
supported by the Australian Research Council, grant DP1094072.


\pagestyle{myheadings} \markright{{Finite index subgroups of mapping class groups}  \quad\sc{\today}}

\part{Preliminaries}

\pagestyle{myheadings} \markright{{Finite index subgroups of mapping class groups}  \quad\sc{\today}}

\section{Useful information on permutations\label{11107S1}}

We define the \emph{minimum index} $\mathrm{mi}(G)$ of a nontrivial group $G$
to be the element of the ordered set $2<3<4<\cdots<\infty$ corresponding to
the minimum among the indices of all proper subgroups of $G$ with finite
index, and $\infty$ when $G$ has no proper finite index subgroup (such a $G$
has been called \emph{counter-finite}).

Permutations are important to our investigation because of the well-known
relationship between index $m$ subgroups and maps to $\mathfrak{S}_{m}$. A
homomorphism $\varphi:G\rightarrow{\mathfrak{S}}_{m}$ is called
\emph{transitive} if its image acts transitively on $\{1,\dots,m\}$. If
$\varphi:G\rightarrow{\mathfrak{S}}_{m}$ is transitive, then $\mathrm{Stab}%
_{\varphi}(1)=\{\gamma\in G\mid\varphi(\gamma)(1)=1\}$ is a subgroup of $G$ of
index $m$. Conversely, if $H$ is a subgroup of $G$ of index $m$, then there
exists a transitive homomorphism $\varphi:G\rightarrow{\mathfrak{S}}_{m}$ such
that $H=\mathrm{Stab}_{\varphi}(1)$ (take the action of $G$ on the right
cosets of $H$). It follows that $\mathrm{mi}(G)$ is also the smallest $m\geq2$
such that there exists a nontrivial homomorphism $\varphi:G\rightarrow
{\mathfrak{S}}_{m}$. Of course, if such a homomorphism does not exist for any
$m$, then $\mathrm{mi}(G)=\infty$.

The minimum index has the property that, if $G\twoheadrightarrow H$ is an
epimorphism, then $\mathrm{mi}(G)\leq\mathrm{mi}(H)$. (Indeed, the pre-image
of an index $m$ subgroup under an epimorphism is an index $m$ subgroup.) From
the definition of $\theta_{g,n}:{\mathcal{M}}_{g,n}\rightarrow\mathrm{Sp}%
_{2g}({\mathbb{F}}_{2})$ as the composition of epimorphisms%
\[
{\mathcal{M}}_{g,n}\twoheadrightarrow{\mathcal{M}}_{g,0}\twoheadrightarrow
\mathrm{Sp}_{2g}({\mathbb{Z}})\twoheadrightarrow\mathrm{Sp}_{2g}({\mathbb{F}%
}_{2}),
\]
we therefore have
\[
\mathrm{mi}({\mathcal{M}}_{g,n})\leq\mathrm{mi}({\mathcal{M}}_{g,0}%
)\leq\mathrm{mi}(\mathrm{Sp}_{2g}({\mathbb{Z}}))\leq\mathrm{mi}(\mathrm{Sp}%
_{2g}({\mathbb{F}}_{2}))\,.
\]

Since many generators of mapping class groups commute with each other, we need
some preliminary results that discuss aspects of commuting permutations. In
the following, $C_{k}$ denotes the cyclic group of order $k$, and an orbit of
cardinality $k$ under the action of a permutation or a permutation group is
called a $k$\emph{-orbit}.

\begin{lemma}
\label{11107L1.1}Let $u\in\mathfrak{S}_{m}$ have cycle type $(1)^{\ell_{1}%
}(2)^{\ell_{2}}\cdots(m)^{\ell_{m}}$, and let
\[
I(u)=\{k\in\{1,2,\dots,m\}\mid\ell_{k}>0\}\,,
\]
so that $\sum_{k\in I(u)}k\ell_{k}=m$.

\begin{enumerate}
\item The centralizer $C_{{\mathfrak{S}}_{m}}(u)$ is isomorphic to
\[
\prod_{k\in I(u)}((C_{k})^{\ell_{k}}\rtimes{\mathfrak{S}}_{\ell_{k}}%
)=\prod_{k\in I(u)}(C_{k}\wr{\mathfrak{S}}_{\ell_{k}})\,.
\]

\item If $P\leq C_{{\mathfrak{S}}_{m}}(u)$ is nonabelian, then, for some $k\in
I(u)$, $\ell_{k}\geq\mathrm{mi}(P)$.

\item If further $P$ is perfect, then $P$ is isomorphic to a subgroup of
\[
\prod_{k\in I(u)}((C_{k})^{\ell_{k}}\rtimes{\mathfrak{A}}_{\ell_{k}}%
)=\prod_{k\in I(u)}(C_{k}\wr{\mathfrak{A}}_{\ell_{k}})\,.
\]
If $I(P)=\{k\in I(u)\mid P\text{ acts nontrivially on the union of the
}k\text{-orbits of }u\}$, then the following numerical constraints hold:

\begin{enumerate}
\item $5\leq\mathrm{mi}(P)\leq\ell_{k}$ whenever $k\in I(P)$;

\item $m\geq\mathrm{mi}(P)(\sum_{k\in I(P)}k)$; and

\item if $m<5\,\mathrm{mi}(P)$, then either $I(P)=\{4\}$ or $I(P)-\{1\}=\{2\}$
or $\{3\}$.
\end{enumerate}
\end{enumerate}
\end{lemma}

\noindent\textbf{Proof.} (1) is left to the reader. Let $P_{k}$ denote the
projection of $P$ to the component $C_{k}\wr{\mathfrak{S}}_{\ell_{k}}$ of
$C_{{\mathfrak{S}}_{m}}(u)$. Since $P$ is contained in $\prod_{k\in I(u)}%
P_{k}$ and $P$ is nonabelian, there exists $k\in I(u)$ such that $P_{k}$ is
nonabelian, and thus $\ell_{k}\geq\mathrm{mi}(P_{k})\geq\mathrm{mi}(P)$.
Finally, (3) results from the fact that any nontrivial perfect subgroup of
$(C_{k})^{\ell_{k}}\rtimes{\mathfrak{S}}_{\ell_{k}}$ (whose support is the
union of the $k$-orbits of $u$) must have nontrivial perfect image in
${\mathfrak{S}}_{\ell_{k}}$. Then $\ell_{k}\geq5$ and the image lies in the
maximum perfect subgroup ${\mathfrak{A}}_{\ell_{k}}$ of ${\mathfrak{S}}%
_{\ell_{k}}$. Then (a), (b) and (c) follow readily. \qed

\bigskip\noindent For $u\in{\mathfrak{S}}_{m}$ we write $\{1,\dots
,m\}=F(u)\sqcup S(u)$ for the partition into the fixed set $F(u)$ and support
$S(u)$ of the permutation $u$. Evidently, $u$ restricts to the identity map on
$F(u)$, respectively to a bijection on $S(u)$.

\begin{lemma}
\label{11107L1.2}Let $u,v\in\mathfrak{S}_{m}$ be such that $uvu=vuv$. Then
$\left\vert S(u)\right\vert \leq2\left\vert S(u)\cap S(v)\right\vert $.
\end{lemma}

\noindent\textbf{Proof.} If $u(i)\in F(v)$, then $(vu)(i)=u(i)$. If also $i\in
F(v)$, then
\[
u(u(i))=(uvu)(i)=(vuv)(i)=(vu)(i)=u(i)\,.
\]
Thus,
\[
i\in u^{-1}(F(v))\cap F(v)\quad\Longrightarrow\quad i\in u^{-1}%
(F(u))\,\text{;}%
\]
in other words, $F(v)\cap u(F(v))\subseteq F(u)$. Therefore,
\[
S(u)\cap F(v)\cap u(F(v))=\emptyset\,\text{;}%
\]
and so $S(u)\cap u(F(v))\subseteq S(u)\cap S(v)$. However, $u$ maps $S(u)\cap
F(v)$ bijectively to $S(u)\cap u(F(v))$. The result is now immediate from the
fact that
\[
\left\vert S(u)\right\vert =\left\vert S(u)\cap S(v)\right\vert +\left\vert
S(u)\cap F(v)\right\vert \text{.}\qed
\]

\begin{lemma}
\label{11107L1.3}Let $u,v\in{\mathfrak{S}}_{m}$ be such that $uvu=vuv$. If
$k\in\{2,3\}$ and all nontrivial orbits of $u$ and $\langle u,v\rangle$ have
length $k$, then $u=v$.
\end{lemma}

\noindent\textbf{Proof.} We present the argument for $k=3$; for $k=2$ it is
similar, but simpler. Since $u$ and $v$ are conjugate, they have the same
cycle decomposition type. If there are no nontrivial orbits, then $u=v=1$. So,
let $O$ be an orbit of $u$ of length $3$. Then $O$ must also be a nontrivial
orbit of $\langle u,v\rangle$. Hence, since $v$ is a product of $3$-cycles,
$v$ acts on $O$ as either $1,u$ or $u^{2}$. In each case, the actions of $u$
and $v$ on $O$ commute. Likewise, on each nontrivial orbit of $v$ the actions
of $u$ and $v$ commute. Finally, on $F(u)\cap F(v)$, since $u$ and $v$ both
act as the identity, the actions of $u$ and $v$ again commute. Hence, $uv=vu$.
From $uvu=vuv$, the result follows. \qed

\begin{lemma}
\label{11107L1.4}Let $u,v_{0},v_{1}\in$\textit{${\mathfrak{S}}$}$_{m}$ be such that

\begin{enumerate}
\item[(a)] $uv_{i}u=v_{i}uv_{i}$ for $i=0,1$; and

\item[(b)] $v_{0}$ and $v_{1}$ commute.
\end{enumerate}

If $u$ has order $3$ and all nontrivial orbits of $\langle u,v_{0}\rangle$ and
$\langle u,v_{1}\rangle$ have length $4$, then $v_{0}=v_{1}$.
\end{lemma}

\noindent\textbf{Proof.} If $u$ is trivial then the result is immediate from
(a). Therefore, we can assume that $u$ contains the $3$-cycle $(1\;2\;3)$ and
$\{2,3,4\}$ is a nontrivial orbit of $v_{0}$. Now, in ${\mathfrak{S}}_{4}$
\[
(1\;2\;3)(4\;3\;2)(1\;2\;3)=(4\;3\;2)(1\;2\;3)(4\;3\;2)=(1\;4)(2\;3)\,;
\]
however
\[
(1\;2\;3)(4\;2\;3)(1\;2\;3)\neq(4\;2\;3)(1\;2\;3)(4\;2\;3)\,.
\]
Thus, $v_{0}$ contains the cycle $(4\;3\;2)$ in its cycle decomposition.
Because $v_{1}$ commutes with $v_{0}$,
\begin{equation}
v_{1}(4\;3\;2)v_{1}^{-1}=(v_{1}(4)\;v_{1}(3)\;v_{1}(2))
\label{Commuting3CycleEqn}%
\end{equation}
is a $3$-cycle in the decomposition of $v_{0}$. Now, $\langle u,v_{1}\rangle$
has a $4$-orbit of the form $\{1,2,3,x\}$, whence, from the above argument in
${\mathfrak{S}}_{4}$, $v_{1}$ acts on this orbit as $(3)(x\;2\;1)$,
$(2)(x\;1\;3)$ or $(1)(x\;3\;2)$. However, in the first case, because
$v_{1}(3)=3$, Equality (\ref{Commuting3CycleEqn}) implies that $v_{1}(2)=2$,
contradicting $v_{1}(2)=1$. Similarly, in the second case, $v_{1}(2)=2$
combines with (\ref{Commuting3CycleEqn}) to imply that $v_{1}(3)=3$,
contradicting $v_{1}(1)=3$.

\smallskip This leaves the last case, in which $v_{1}(3)=2$, which from
(\ref{Commuting3CycleEqn}) gives $x=v_{1}(2)=4$. That is, $(4\;3\;2)$ is a
$3$-cycle in the decomposition of $v_{1}$. It follows that the $3$-cycles of
$v_{0}$ and $v_{1}$ coincide, whence $v_{0}=v_{1}$. \qed

\begin{lemma}
\label{11107L1.5}Let $u,v_{0},v_{1},v_{2}\in\mathfrak{S}_{m}$ be such that

\begin{enumerate}
\item[(a)] $uv_{i}u=v_{i}uv_{i}$ for $i=0,1,2$;

\item[(b)] $v_{0}$, $v_{1}$ and $v_{2}$ commute pairwise; and

\item[(c)] whenever $\{i,j,k\}=\{0,1,2\}$, there is an isomorphism
$\gamma_{j,k}$ from $\langle u,v_{i},v_{j}\rangle$ to $\langle u,v_{i}%
,v_{k}\rangle$ fixing $u$ and $v_{i}$, and sending $v_{j}$ to $v_{k}$.
\end{enumerate}

If $u$ is a product of $4$-cycles (possibly with fixed points) and all
nontrivial orbits of each $\langle u,v_{i}\rangle$ have length $4$, then
$v_{0}=v_{1}=v_{2}$.
\end{lemma}

\noindent\textbf{Proof.} Since, by (a), $u$ and $v_{i}$ are conjugate, they
have the same cycle decomposition type, which must be a product of $4$-cycles.
If there are no nontrivial orbits, then each $v_{i}=1$.

\smallskip Let $O$ be an orbit of $u$ of length $4$. Then $O$ must also be a
nontrivial orbit of each $\langle u,v_{i}\rangle$. From (a) it follows that
$O$ is also a nontrivial orbit of each $v_{i}$. Likewise, a nontrivial orbit
of any one $v_{i}$ must also be a nontrivial orbit of $u$ and thereby of each
$v_{i}$. Thus, it suffices to restrict attention to the action on each
nontrivial orbit $O$; in effect, $\langle u,v_{0},v_{1},v_{2}\rangle
\leq{\mathfrak{S}}_{4}$. Using (b), from the fact that in ${\mathfrak{S}}_{4}$
commuting $4$-cycles are either the same or mutually inverse, we must have at
least two distinct indices $i,j\in\{0,1,2\}$ such that $v_{i}=v_{j}$. Since by
(c) $\gamma_{j,k}$ both fixes $v_{i}$ and sends $v_{j}$ to $v_{k}$, we
conclude that $v_{i}=v_{j}=v_{k}$. \qed


\section{Presentation for the mapping class group\label{11107S2}}

Throughout the paper we denote by $T_{b}$ the Dehn twist about a simple closed
curve $b$. We fix a representation of $\Sigma_{g,1}$ as well as the simple
closed curves $a_{0}, a_{1}, \dots, a_{2g+1}$ illustrated in Figure~2.1, and
we set $T_{i} = T_{a_{i}}$ for all $0 \le i \le2g+1$. The following result is
shown in \cite{Matsu1}.

\begin{figure}[tbh]
\bigskip\centerline{
\setlength{\unitlength}{0.5cm}
\begin{picture}(24.5,4)
\put(0,0){\includegraphics[width=12.25cm]{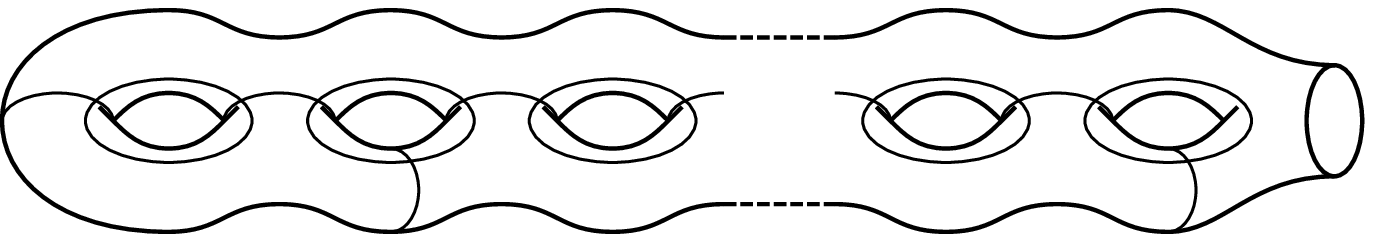}}
\put(1,2.8){\small $a_1$}
\put(2.9,3){\small $a_2$}
\put(5,2.8){\small $a_3$}
\put(6.9,3){\small $a_4$}
\put(7.7,0.7){\small $a_0$}
\put(9,2.8){\small $a_5$}
\put(10.9,3){\small $a_6$}
\put(16.3,3.1){\small $a_{2g-2}$}
\put(18.3,2.9){\small $a_{2g-1}$}
\put(20.8,3){\small $a_{2g}$}
\put(19.7,0.7){\small $a_{2g+1}$}
\end{picture}} \bigskip
\centerline{{\bf Figure 2.1.} Generators for $\mathcal{M}_{g,1}$.}
\bigskip\end{figure}

\begin{theorem}
\label{11107T2.1}\emph{(Matsumoto \cite{Matsu1})}. Let $g\geq2$.

\begin{enumerate}
\renewcommand{\labelenumi}{\alph{enumi})}

\item $\mathcal{M}_{g,1}$ has a presentation with generators $T_{0}%
,T_{1},\dots,T_{2g}$ and relations
\begin{gather*}
T_{i}T_{j}T_{i}=T_{j}T_{i}T_{j}\quad\text{if }a_{i}\text{ and }a_{j}\text{
intersect in a single point,}\\
T_{i}T_{j}=T_{j}T_{i}\quad\text{if }a_{i}\cap a_{j}=\emptyset\,,\\
(T_{2}T_{3}T_{4}T_{0})^{10}=(T_{1}T_{2}T_{3}T_{4}T_{0})^{6}\,,\\
(T_{2}T_{3}T_{4}T_{5}T_{6}T_{0})^{12}=(T_{1}T_{2}T_{3}T_{4}T_{5}T_{6}%
T_{0})^{9}\quad\text{if }g\geq3\,.
\end{gather*}

\item $\mathcal{M}_{g,0}$ is the quotient of $\mathcal{M}_{g,1}$ by the
additional relation
\[
T_{1}^{2g-2} = (T_{0} T_{3} T_{4}\cdots T_{2g-1})^{4g-4} .
\]

\end{enumerate}

\vspace{-0.25in}\qed
\end{theorem}

\medskip We call $k$ simple closed curves $b_{1},\dots,b_{k}$ in $\Sigma
_{g,n}$ a $k$-\emph{chain} if their intersection numbers satisfy
$i(b_{i},b_{j})=1$ if $|i-j|=1$ and $i(b_{i},b_{j})=0$ otherwise. Such a
$k$-chain is called \emph{nonseparating} if the complement of $b_{1}\cup
\cdots\cup b_{k}$ in $\Sigma_{g,n}$ is connected. Note that, if $(b_{1}%
,\dots,b_{k})$ is a nonseparating $k$-chain, then each $b_{i}$ is
nonseparating, too. The next lemma follows from the rigidity of closed Riemann
surfaces (see, for example, \cite{FarMar1} Sections 1.3 and 2.3).

\begin{lemma}
\label{11107L2.2}Let $(b_{1},\dots,b_{k})$ and $(b_{1}^{\prime},\dots
,b_{k}^{\prime})$ be two nonseparating $k$-chains in $\Sigma_{g,n}$. Then
there exists $\alpha\in\mathcal{M}_{g,n}$ such that $\alpha(b_{i}%
)=b_{i}^{\prime}$ for all $1\leq i\leq k$. In consequence, this $\alpha$
satisfies $\alpha T_{b_{i}}\alpha^{-1}=T_{b_{i}^{\prime}}$ for all $1\leq
i\leq k$. \qed

\end{lemma}

\medskip For studying the mapping class group ${\mathcal{M}}_{g,n}$ with
$n\geq2$ we use the following convention: $a_{0},a_{2},\dots,a_{2g},a_{2g+1}$
and $b_{1},\dots,b_{n}$ are the simple closed curves illustrated in Figure
2.2, $T_{i}=T_{a_{i}}$ for all $0\leq i\leq2g+1$, $i\neq1$, and $T_{j}%
^{\prime}=T_{b_{j}}$ for all $1\leq j\leq n$. In order to unify statements for
$n=1$ and $n>1$, we make the further convention that, when $n=1$, $b_{1}$ in
Figure 2.2 coincides with $a_{1}$ in Figure 2.1. Then for $n=1$ the element
$T_{1}^{\prime}$ is simply $T_{1}$.

\begin{figure}[tbh]
\bigskip\centerline{
\setlength{\unitlength}{0.5cm}
\begin{picture}(26,12)
\put(0,0){\includegraphics[width=13cm]{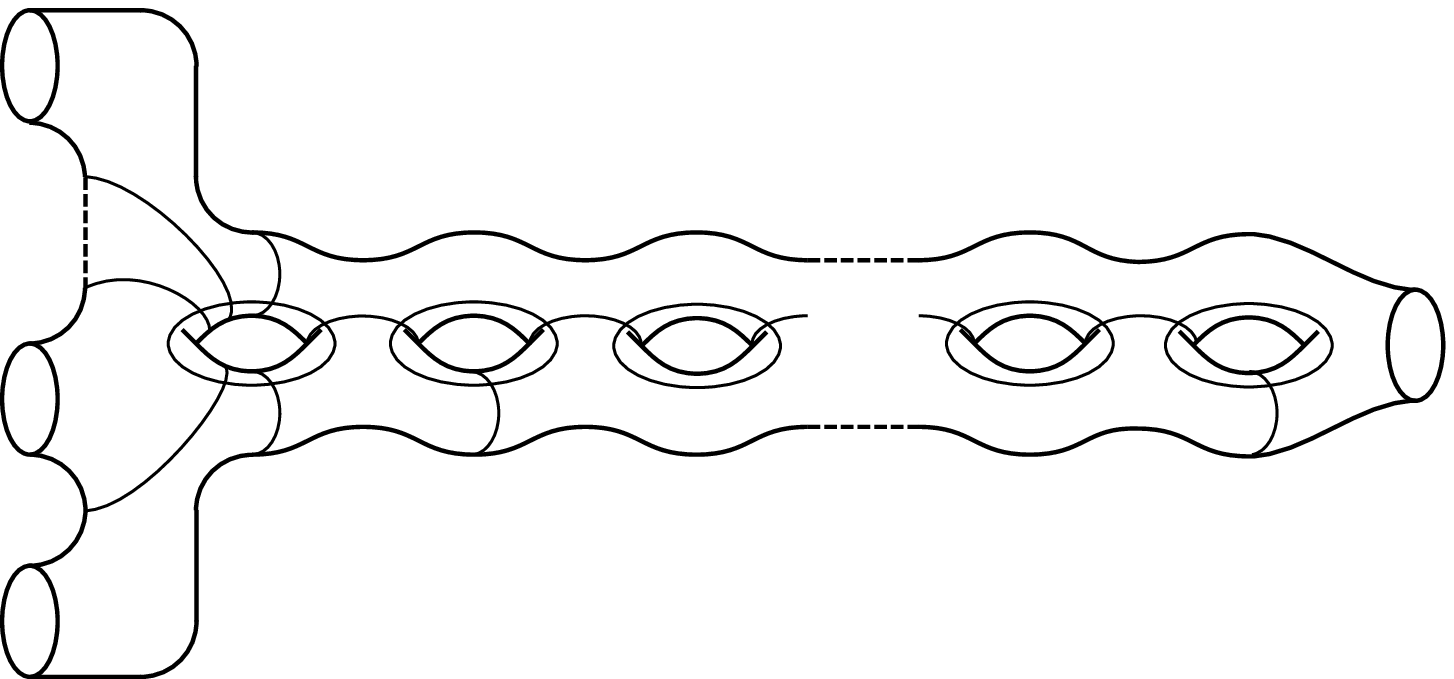}}
\put(5.1,4.6){\small $b_1$}
\put(2.2,4){\small $b_2$}
\put(1.9,7.4){\small $b_3$}
\put(2,9){\small $b_{n-1}$}
\put(5.1,7){\small $b_n$}
\put(2.2,5.7){\small $a_2$}
\put(6.1,6.7){\small $a_3$}
\put(8.1,7){\small $a_4$}
\put(9,4.6){\small $a_0$}
\put(10,6.8){\small $a_5$}
\put(12,6.9){\small $a_6$}
\put(18,7.1){\small $a_{2g-2}$}
\put(20,6.8){\small $a_{2g-1}$}
\put(22.5,7){\small $a_{2g}$}
\put(21.1,4.6){\small $a_{2g+1}$}
\end{picture}} \bigskip
\centerline{{\bf Figure 2.2.} Generators for $\MM_{g,n}$.} \bigskip
\end{figure}

\bigskip\noindent The following is well-known. (It can be found for instance
in \cite[Prop. 2.10 and Thm 3.1]{LabPar1}.

\begin{proposition}
\label{11107P2.3}Let $g\geq2$ and $n\geq1$. Then \textit{${\mathcal{M}}$%
}$_{g,n}$ is generated by $T_{0},T_{2},\dots,T_{2g},T_{1}^{\prime},\dots
,T_{n}^{\prime}$. \qed

\end{proposition}

\medskip\textbf{Observation.} Let $g\geq3$ and $n\geq1$. There is an injective
homomorphism ${\mathcal{M}}_{g-1,n}\rightarrow{\mathcal{M}}_{g,n}$ which sends
$T_{j}^{\prime}$ to $T_{j}^{\prime}$ and $T_{i}$ to $T_{i}$ for all
$j\in\{1,\dots,n\}$ and $i\in\{0,2,\dots,2g-2\}$. It is easily seen that this
homomorphism is induced by some embedding of $\Sigma_{g-1,n}$ into
$\Sigma_{g,n}$. From now on we will assume ${\mathcal{M}}_{g-1,n}$ to be
embedded into ${\mathcal{M}}_{g,n}$ via this homomorphism. Note that such a
homomorphism does not exist for $n=0$.


\section{The genus $2$ case\label{11107S3}}

In this section we describe all the subgroups of ${\mathcal{M}}_{2,n}$ of
index at most $N_{2}^{+} =10$ up to conjugation. We will see in particular
that the genus $g=2$ case is different from the genus $g \ge3$ case.

\medskip The first difference comes from the fact that the abelianization of
${\mathcal{M}}_{2,n}$ is nontrivial, while the group ${\mathcal{M}}_{g,n}$ is
perfect if $g\geq3$ (see \cite{Powel1}, \cite{Korkm1}). More precisely, the
abelianization of ${\mathcal{M}}_{2,n}$ is ${\mathbb{Z}}/10{\mathbb{Z}}$
\cite{Mumfo1}. So, if $\mathrm{ab}:{\mathcal{M}}_{2,n}\rightarrow{\mathbb{Z}%
}/10{\mathbb{Z}}$ denotes the projection of ${\mathcal{M}}_{2,n}$ onto its
abelianization, then $\mathrm{ab}^{-1}({\mathbb{Z}}/5{\mathbb{Z}})$ is a
subgroup of index $2$, $\mathrm{ab}^{-1}({\mathbb{Z}}/2{\mathbb{Z}})$ is a
subgroup of index $5$, and $\mathrm{Ker}(\mathrm{ab})$ is a subgroup of index
$10$. Note that all these subgroups contain the commutator subgroup
${\mathcal{M}}_{2,n}^{\prime}$. In particular, since $\mathrm{Sp}%
_{4}({\mathbb{F}}_{2})^{\prime}$ has index 2 in $\mathrm{Sp}_{4}({\mathbb{F}%
}_{2})$, neither $\mathrm{ab}^{-1}({\mathbb{Z}}/2{\mathbb{Z}})$ nor
$\mathrm{Ker}(\mathrm{ab})$ is a pre-image under $\theta_{2,n}$ of a subgroup
of $\mathrm{Sp}_{4}({\mathbb{F}}_{2})$.

\medskip The second difference comes from the fact that $\mathrm{Sp}%
_{4}({\mathbb{F}}_{2})={\mathfrak{S}}_{6}$ \cite{DicksonTAMS1908} contains
more than two subgroups of index at most $10$, up to conjugation. In addition
to $O_{4}^{-}({\mathbb{F}}_{2})={\mathfrak{S}}_{5}$ \cite{DicksonTAMS1908} and
$O_{4}^{+}({\mathbb{F}}_{2})$, of indices $6$ and $10$ respectively, it
contains the alternating group ${\mathfrak{A}}_{6}$ of index $2$, and another
subgroup of index $6$ which can be described as follows. The group
$\mathrm{Sp}_{4}({\mathbb{F}}_{2})={\mathfrak{S}}_{6}$ has a non-inner
automorphism $\alpha$ defined by
\[
\alpha:\left\{
\begin{array}
[c]{l}%
(1\;2)\mapsto(1\;2)(3\;5)(4\;6)\\
(2\;3)\mapsto(1\;3)(2\;4)(5\;6)\\
(3\;4)\mapsto(1\;2)(3\;6)(4\;5)\\
(4\;5)\mapsto(1\;3)(2\;5)(4\;6)\\
(5\;6)\mapsto(1\;2)(3\;4)(5\;6)
\end{array}
\right.
\]
(Since $\mathrm{Out}({\mathfrak{S}}_{6})$ has order $2$, $\alpha$ is
essentially unique.) It turns out that $\alpha(O_{4}^{-}({\mathbb{F}}_{2}))$
is a subgroup of $\mathrm{Sp}_{4}({\mathbb{F}}_{2})$ of index $6$ which is not
conjugate to $O_{4}^{-}({\mathbb{F}}_{2})$. On the other hand, $\alpha
(O_{4}^{+}({\mathbb{F}}_{2}))$ is conjugate to $O_{4}^{+}({\mathbb{F}}_{2})$
(so there are four subgroups and not five), and $\mathrm{ab}^{-1}({\mathbb{Z}%
}/5{\mathbb{Z}})=\theta_{2}^{-1}({\mathfrak{A}}_{6})$. Note also that
$\mathrm{Out}({\mathfrak{S}}_{n})$ is trivial if $n\neq6$.

\medskip So, we have the following subgroups of index at most $10$ in
${\mathcal{M}}_{2,n}$ up to conjugation:

\begin{itemize}
\item $\theta_{2,n}^{-1}({\mathfrak{A}}_{6})=\mathrm{ab}^{-1}({\mathbb{Z}%
}/5{\mathbb{Z}})$ of index $2$,

\item $\mathrm{ab}^{-1}({\mathbb{Z}}/2{\mathbb{Z}})$ of index $5$,

\item $\theta_{2,n}^{-1}(O_{4}^{-}({\mathbb{F}}_{2}))$ of index $6$,

\item $\theta_{2,n}^{-1}(\alpha(O_{4}^{-}({\mathbb{F}}_{2})))$ of index $6$,

\item $\mathrm{Ker}(\mathrm{ab})$ of index $10$, and

\item $\theta_{2,n}^{-1}(O_{4}^{+}({\mathbb{F}}_{2}))$ of index $10$.
\end{itemize}

We show that these are all the proper subgroups of index at most $10$ in
${\mathcal{M}}_{2,n}$ up to conjugation.

\medskip Let $n\geq1$, $m\geq1$ and $w\in{\mathfrak{S}}_{m}$. If $w^{10}=1$,
then there is a permutation representation $\mathrm{cycl}_{w}:{\mathcal{M}%
}_{2,n}\rightarrow{\mathfrak{S}}_{m}$ which sends $T_{i}$ and $T_{j}^{\prime}$
to $w$ for all $i\in\{0,2,3,4\}$ and all $j\in\{1,\dots,n\}$. Such a
representation is called a \emph{cyclic representation} of ${\mathcal{M}%
}_{2,n}$. It is transitive if and only if $w$ is a cycle of length $m$ and
$m\in\{1,2,5,10\}$.

\begin{lemma}
\label{11107L3.1}\label{11107L3.2} For $n\geq1$, there are exactly three
conjugacy classes of non-cyclic transitive permutation representations of
$\mathcal{M}_{2,n}$ of degree at most $10$, namely, two conjugacy classes of
degree~$6$, and a unique conjugacy class of degree $10$. More specifically, up
to conjugacy in $\mathfrak{S}_{6}$, the two transitive permutation
representations $\mathcal{M}_{2,n}\rightarrow\mathfrak{S}_{6}$ are given by
\begin{align*}
\phi_{2,n}^{-}:  &  \left\{
\begin{array}
[c]{ll}%
T_{j}^{\prime}\mapsto(1\;2) & \text{ for }1\leq j\leq n\\
T_{2}\mapsto(2\;3) & \\
T_{3}\mapsto(3\;4) & \\
T_{4}\mapsto(4\;5) & \\
T_{0}\mapsto(5\;6) &
\end{array}
\right. \\[1ex]
\phi^{\alpha}_{n}:  &  \left\{
\begin{array}
[c]{ll}%
T_{j}^{\prime}\mapsto(1\;2)(3\;5)(4\;6) & \text{ for }1\leq j\leq n\\
T_{2}\mapsto(1\;3)(2\;4)(5\;6) & \\
T_{3}\mapsto(1\;2)(3\;6)(4\;5) & \\
T_{4}\mapsto(1\;3)(2\;5)(4\;6) & \\
T_{0}\mapsto(1\;2)(3\;4)(5\;6) &
\end{array}
\right.
\end{align*}
and, up to conjugacy in $\mathfrak{S}_{10}$, the unique permutation
representation $\mathcal{M}_{2,n}\rightarrow\mathfrak{S}_{10}$ is given by
\[
\phi_{2,n}^{+}:\left\{
\begin{array}
[c]{ll}%
T_{j}^{\prime}\mapsto(3\;5)(6\;8)(9\;10) & \text{ for }1\leq j\leq n\\
T_{2}\mapsto(2\;3)(4\;6)(7\;9) & \\
T_{3}\mapsto(1\;2)(6\;10)(8\;9) & \\
T_{4}\mapsto(2\;4)(3\;6)(5\;8) & \\
T_{0}\mapsto(4\;7)(6\;9)(8\;10) &
\end{array}
\right.
\]

\end{lemma}

\noindent\textbf{Proof.} In the case $n=1$, this result can be easily proved
with a direct calculation. Using the presentation of ${\mathcal{M}}_{2,1}$
from Theorem \ref{11107T2.1}, one can use coset enumeration techniques to
perform a systematic search for representatives of the conjugacy classes of
subgroups of ${\mathcal{M}}_{2,1}$ of index at most $K$ for a (small) integer
$K$; see \cite{Sims1}. When $K=10$, this systematic search shows that there
are exactly six conjugacy classes of proper subgroups of ${\mathcal{M}}_{2,1}$
of index at most $10$; the columns of the coset table for the each of the
constructed subgroups yield the images of the generators $T_{0},\dots,T_{4}$
under the corresponding permutation representation. The computation is very
easy and can be performed with any mathematical software such as
\textsc{Magma} or GAP.

Now suppose that $n\geq2$. Let $\varphi:{\mathcal{M}}_{2,n}\rightarrow
{\mathfrak{S}}_{m}$ be a non-cyclic and transitive representation with
$m\leq10$. For $1\leq j\leq n$, we denote by ${\mathcal{M}}^{(j)}$ the
subgroup of ${\mathcal{M}}_{2,n}$ generated by $T_{j}^{\prime},T_{2}%
,T_{3},T_{4},T_{0}$. This group is isomorphic to ${\mathcal{M}}_{2,1}$ via an
isomorphism $\gamma_{j}:{\mathcal{M}}_{2,1}\rightarrow{\mathcal{M}}^{(j)}$
which sends $T_{1}$ to $T_{j}^{\prime}$ and $T_{i}$ to $T_{i}$ for all
$i\in\{2,3,4,0\}$. We denote by $\varphi_{j}:{\mathcal{M}}_{2,1}%
\rightarrow{\mathfrak{S}}_{m}$ the composition of $\gamma_{j}$ with $\varphi$.
Observe that

\begin{itemize}
\item[($\ast$)] $\varphi_{j}(T_{i})=\varphi(T_{i})$ for all $j\in
\{1,\dots,n\}$ and $i\in\{0,2,3,4\}$.
\end{itemize}

By the case $n=1$, for each $j\in\{1,\dots,n\}$ there is a decomposition
$\{1,\dots,m\}=S_{j}^{(1)}\sqcup S_{j}^{(2)}\sqcup S_{j}^{(3)}$ as follows.
Each $S_{j}^{(k)}$ is invariant under the action of $\varphi_{j}({\mathcal{M}%
}_{2,1})$; either $S_{j}^{(1)}=\emptyset$, or $\varphi_{j}$ restricted to
$S_{j}^{(1)}$ is equivalent to an element of $\{\phi_{2,1}^{-},\phi
_{1}^{\alpha},\phi_{2,1}^{+}\}$; if $S_{j}^{(2)}$ is nonempty, then there
exists $w_{j}\in{\mathfrak{S}}_{m}-\{1\}$ such that $S(w_{j})=S_{j}^{(2)}$ and
the restriction of $\varphi_{j}$ to $S_{j}^{(2)}$ is $\mathrm{cycl}_{w_{j}}$;
and $\varphi_{j}({\mathcal{M}}_{2,1})$ acts trivially on $S_{j}^{(3)}$.

Let $\phi\in\{\phi_{2,1}^{-}, \phi_{1}^{\alpha}, \phi_{2,1}^{+}\}$ and set
$N=N_{\phi}=6$ if $\phi=\phi_{2,1}^{-}$ or $\phi_{1}^{\alpha}$, and
$N=N_{\phi}=10$ if $\phi=\phi_{2,1}^{+}$. The following claims are readily
verified from the description of $\phi$, using either GAP or \textsc{Magma}
where necessary:

\begin{enumerate}
\item $\phi(T_{2})$ and $\phi(T_{3})$ have no common cycle in their decompositions;

\item $S(\{\phi(T_{i})\mid i=1,3,4,0\})=\{1,\dots,N\}$, where, for
$X\subseteq{\mathfrak{S}}_{N}$, $S(X)$ denotes $\cup_{w\in X}S(w)$;

\item the (simultaneous) centralizer of $\{\phi(T_{i})\mid i=1,3,4,0\}$ in
${\mathfrak{S}}_{N}$ is $\{1,\phi(T_{1})\}$;

\item the support of each cycle in the decomposition of $\phi(T_{1})$
intersects $S(\{\phi(T_{i})\mid i=0,2,3,4\})$ nontrivially.
\end{enumerate}

We first show that each $S_{j}^{(2)}=\emptyset$. Whenever $S_{j}^{(2)}%
\neq\emptyset$, then by (1) we get that $w_{j}$ is the product of the common
nontrivial cycles of $\varphi_{j}(T_{2})$ and $\varphi_{j}(T_{3})$. However,
by ($\ast$), these common cycles are independent of choice of $j$. It follows
that for all $j\in\{1,\dots,n\}$, $S_{j}^{(2)}\neq\emptyset$, $w_{j}=w_{1}$
and $S_{j}^{(2)}=S(w_{j})=S(w_{1})=S_{1}^{(2)}$. However, $\varphi$ is
transitive and non-cyclic. Hence, $S_{j}^{(2)}=\emptyset$ for all
$j\in\{1,\dots,n\}$.

It now follows that each $S_{j}^{(1)}\neq\emptyset$, and therefore, the
restriction of each $\varphi_{j}$ to $S_{j}^{(1)}$ is equivalent to an element
of $\{\phi_{2,1}^{-},\phi_{1}^{\alpha},\phi_{2,1}^{+}\}$. Since we know by
($\ast$) that for each $i\neq1$ the $\varphi_{j}(T_{i})$ agree, it remains to
show that all the $\varphi_{j}(T_{1})$ also coincide, and that each
$S_{j}^{(3)}$ is empty.

Without loss of generality we can assume that $S_{1}^{(1)}=\{1,\dots,N\}$ and
the restriction of $\varphi_{1}$ to $\{1,\dots,N\}$ is an element $\phi$ in
$\{\phi_{2,1}^{-},\phi_{1}^{\alpha},\phi_{2,1}^{+}\}$, where $N=N_{\phi}$. Let
$j\in\{1,\dots,n\}$. Since $T_{j}^{\prime}$ commutes with $T_{1}^{\prime
},T_{3},T_{4},T_{0}$, the permutation $\varphi_{j}(T_{1})$ belongs to the
centralizer of $\{\varphi_{1}(T_{i})\mid i=1,3,4,0\}$. Combining (2) and (3)
we get that this centralizer is $\{1,\phi(T_{1})\}\times{\mathfrak{S}}_{m-N}$.
On the other hand, by (4), the support of each cycle of $\varphi_{j}(T_{1})$
intersects $S(\{\varphi_{j}(T_{i})\mid i=0,2,3,4\})$ nontrivially, and, by
($\ast$),
\[
S(\{\varphi_{j}(T_{i})\mid i=0,2,3,4\})=S(\{\varphi_{1}(T_{i})\mid
i=0,2,3,4\})\subseteq\{1,\dots,N\}.
\]
This implies that $\varphi_{j}(T_{1})=\phi(T_{1})$ and, therefore, each
$\varphi_{j}(T_{1})$ is the same. Then, finally, by transitivity, it must be
that $m=N$, and the proof is complete.\qed

\begin{proposition}
\label{11107P3.3}Let $n\geq0$. Then $\mathcal{M}_{2,n}$ has precisely six
proper subgroups of index at most $10$ up to conjugation, namely,
$\mathrm{ab}^{-1}(\mathbb{Z}/5\mathbb{Z})=\theta_{2,n}^{-1}(\mathfrak{A}_{6})$
of index $2$, $\mathrm{ab}^{-1}(\mathbb{Z}/2\mathbb{Z})$ of index $5$,
$\theta_{2,n}^{-1}(O_{4}^{-}(\mathbb{F}_{2}))$ of index $6$, $\theta
_{2,n}^{-1}(\alpha(O_{4}^{-}(\mathbb{F}_{2})))$ of index $6$, $\mathrm{Ker}%
(\mathrm{ab})$ of index $10$, and $\theta_{2,n}^{-1}(O_{4}^{+}(\mathbb{F}%
_{2}))$ of index $10$.
\end{proposition}

\noindent\textbf{Proof.} For $n\ge1$, the claim is a direct consequence of
Lemma \ref{11107L3.1} together with the description of the subgroups of
${\mathcal{M}}_{2,n}$ given in the beginning of the section. By
Theorem~\ref{11107T2.1}, ${\mathcal{M}}_{2,0}$ is a quotient of ${\mathcal{M}%
}_{2,1}$ by one additional relation, so for the case $n=0$ it is sufficient to
check that the representations of ${\mathcal{M}}_{2,1}$ given in Lemma
\ref{11107L3.1} satisfy the additional relation of ${\mathcal{M}}_{2,0}$; this
is the case, as can easily be verified. \qed


\section{The genus $3$ case\label{11107S4}}

In this section we calculate the subgroups of index at most $36$ in
${\mathcal{M}}_{3,n}$ up to conjugation. (Note that $N_{3}^{-}=28$ and
$N_{3}^{+}=36$.) We argue in the same manner as for the case of genus~$2$
surfaces (see Section \ref{11107S3}), with direct calculations often made with
computers. However, we should point out here that, in this case, the
computations are far from being elementary, and we often approach the limit of
what can currently be done with computers (especially in the proof of
Lemma~\ref{11107L4.1}). Recall also that the case of surfaces of genus $g=3$
will be the first step in the induction to prove Theorem~\ref{11107T0.2}.

\begin{lemma}
\label{11107L4.1}\label{11107L4.2}For $n\geq1$, there are exactly two
conjugacy classes of nontrivial transitive permutation representations of
$\mathcal{M}_{3,n}$ of degree at most $36$, namely a unique conjugacy class of
degree~$28$ and a unique conjugacy class of degree $36$. More specifically, up
to conjugacy in $\mathfrak{S}_{28}$, the unique permutation representation
$\mathcal{M}_{3,n}\rightarrow\mathfrak{S}_{28}$ is given by
\[
\phi_{3,n}^{-}:\left\{
\begin{array}
[c]{l}%
T_{j}^{\prime}\mapsto(14\;18)(16\;21)(17\;22)(19\;23)(24\;26)(27\;28)\text{
for }1\leq j\leq n\\
T_{0}\mapsto(1\;3)(2\;5)(4\;8)(20\;25)(24\;28)(26\;27)\\
T_{2}\mapsto(10\;14)(12\;16)(13\;17)(15\;19)(20\;24)(25\;28)\\
T_{3}\mapsto(6\;10)(7\;12)(9\;13)(11\;15)(24\;27)(26\;28)\\
T_{4}\mapsto(3\;6)(4\;7)(5\;9)(15\;20)(19\;24)(23\;26)\\
T_{5}\mapsto(2\;4)(5\;8)(6\;11)(10\;15)(14\;19)(18\;23)\\
T_{6}\mapsto(1\;2)(3\;5)(6\;9)(10\;13)(14\;17)(18\;22)
\end{array}
\right.
\]
and, up to conjugacy in $\mathfrak{S}_{36}$, the unique permutation
representation $\mathcal{M}_{3,n}\rightarrow\mathfrak{S}_{36}$ is given by
\[
\phi_{3,n}^{+}:\left\{
\begin{array}
[c]{l}%
T_{j}^{\prime}\mapsto
(6\;9)(10\;13)(14\;18)(15\;19)(20\;22)(21\;25)(26\;28)(27\;30)(31\;32)(34\;35)\\
\hspace{11cm}\text{for }1\leq j\leq n\\
T_{0}\mapsto
(1\;2)(11\;17)(14\;22)(16\;24)(18\;20)(21\;28)(23\;29)(25\;26)(27\;32)(30\;31)\\
T_{2}\mapsto
(4\;6)(7\;10)(11\;14)(12\;15)(16\;21)(17\;22)(23\;27)(24\;28)(29\;32)(33\;35)\\
T_{3}\mapsto
(3\;4)(5\;7)(8\;12)(14\;20)(18\;22)(21\;26)(25\;28)(27\;31)(30\;32)(33\;36)\\
T_{4}\mapsto
(2\;3)(7\;11)(10\;14)(12\;16)(13\;18)(15\;21)(19\;25)(29\;33)(31\;34)(32\;35)\\
T_{5}\mapsto
(3\;5)(4\;7)(6\;10)(9\;13)(16\;23)(21\;27)(24\;29)(25\;30)(26\;31)(28\;32)\\
T_{6}\mapsto
(5\;8)(7\;12)(10\;15)(11\;16)(13\;19)(14\;21)(17\;24)(18\;25)(20\;26)(22\;28)
\end{array}
\right.
\]
In particular, $\mathrm{mi}(\mathcal{M}_{3,n})=28$.
\end{lemma}

\noindent\textbf{Proof.} In the case $n=1$, the result is shown by a direct
computation. Using the presentation of ${\mathcal{M}}_{3,1}$ from Theorem
\ref{11107T2.1}, one can use coset enumeration techniques to perform a
systematic search for representatives of the conjugacy classes of subgroups of
${\mathcal{M}}_{3,1}$ of index at most $K$ for a (small) integer $K$; see
\cite{Sims1}. When $K=36$, this systematic search shows that there are exactly
two conjugacy classes of proper subgroups of ${\mathcal{M}}_{3,1}$ of index at
most $36$: exactly one conjugacy class of subgroups of index $28$ and exactly
one conjugacy class of subgroups of index~$36$. The columns of the coset table
for the each of the constructed subgroups yield the images of the generators
$T_{0},\dots,T_{6}$ under the corresponding permutation representations
$\phi_{3,1}^{-}$ of degree~$28$ and $\phi_{3,1}^{+}$ of degree $36$.

\noindent We used the implementation of the low index subgroup search provided
in \textsc{Magma} \cite{BoCaPl1}, filling the coset table in column major
order. We ran a development version of \textsc{Magma} V2.15. The computation
took approximately 47.5 hours on a GNU\thinspace/\thinspace Linux system with
an Intel E8400 64-bit CPU (core:~3\thinspace GHz, FSB:~1333\thinspace MHz) and
a main memory bandwidth of 6.5\thinspace GB/s (X38 chipset, dual channel DDR2
RAM, memory bus: 1066\thinspace MHz). We remark that the use of column major
order is crucial for the running time; tests for indices between $10$ and $15$
suggest a speed-up by a factor between $10^{3}$ and $10^{4}$ compared to row
major order.

Now suppose that $n\geq2$. Let $\varphi:{\mathcal{M}}_{3,n}\rightarrow
{\mathfrak{S}}_{m}$ be a transitive representation with $m\leq36$. For $1\leq
j\leq n$, we denote by ${\mathcal{M}}^{(j)}$ the subgroup of ${\mathcal{M}%
}_{3,n}$ generated by $T_{j}^{\prime},T_{2}, \dots, T_{6}, T_{0}$. This group
is isomorphic to ${\mathcal{M}}_{3,1}$ via an isomorphism $\gamma
_{j}:{\mathcal{M}}_{3,1}\rightarrow{\mathcal{M}}^{(j)}$ which sends $T_{1}$ to
$T_{j}^{\prime}$ and $T_{i}$ to $T_{i}$ for all $i\in\{2,\dots, 6,0\}$. We
denote by $\varphi_{j}:{\mathcal{M}}_{3,1}\rightarrow{\mathfrak{S}}_{m}$ the
composition of $\gamma_{j}$ with $\varphi$. Observe that

\begin{itemize}
\item[($\ast$)] $\varphi_{j}(T_{i})=\varphi(T_{i})$ for all $j\in
\{1,\dots,n\}$ and $i\in\{0,2, \dots, 6\}$.
\end{itemize}

By the case $n=1$, for each $j\in\{1,\dots,n\}$ there is a decomposition
$\{1,\dots,m\}=S_{j}^{(1)}\sqcup S_{j}^{(2)}$ as follows. Each $S_{j}^{(k)}$
is invariant under the action of $\varphi_{j}({\mathcal{M}}_{3,1})$; the
restriction of $\varphi_{j}$ to $S_{j}^{(1)}$ is equivalent to an element of
$\{\phi_{3,1}^{-},\phi_{3,1}^{+}\}$; and $\varphi_{j}({\mathcal{M}}_{3,1})$
acts trivially on $S_{j}^{(2)}$. (Note here a difference from the $g=2$ case:
because ${\mathcal{M}}_{3,1}$ is a perfect group \cite{Powel1} it has no
cyclic representations.)

Let $\phi\in\{\phi_{3,1}^{-},\phi_{3,1}^{+}\}$ and set $N=N_{\phi}=28$ if
$\phi=\phi_{3,1}^{-}$, and $N=N_{\phi}=36$ if $\phi=\phi_{3,1}^{+}$. Again,
the following claims may be confirmed from the description of $\phi$, using
either GAP or \textsc{Magma} where necessary:

\begin{enumerate}
\item $S(\{\phi(T_{i})\mid i=1,3,4,5,6,0\})=\{1,\dots,N\}$, where, for
$X\subseteq{\mathfrak{S}}_{N}$, $S(X)$ denotes $\cup_{w\in X}S(w)$;

\item the unique element of the centralizer of $\{\phi(T_{i})\mid
i=1,3,4,5,6,0\}$ in ${\mathfrak{S}}_{N}$ having the same cycle decomposition
type as $\phi(T_{1})$ is $\phi(T_{1})$;

\item the support of each nontrivial cycle in the decomposition of $\phi
(T_{1})$ intersects $S(\{\phi(T_{i})\mid i=0,2,3,4,5,6\})$ nontrivially.
\end{enumerate}

Without loss of generality we can assume that $S_{1}^{(1)}=\{1,\dots,N\}$ and
the restriction of $\varphi_{1}$ to $\{1,\dots,N\}$ is an element $\phi$ in
$\{\phi_{3,1}^{-},\phi_{3,1}^{+}\}$, where $N=N_{\phi}$. Let $j\in
\{1,\dots,n\}$. Since $T_{j}^{\prime}$ commutes with $T_{1}^{\prime}%
,T_{3},T_{4},T_{5},T_{6},T_{0}$, the permutation $\varphi_{j}(T_{1})$ belongs
to the centralizer of $\{\varphi_{1}(T_{i})\mid i=1,3,4,5,6,0\}$. By (1), this
centralizer is $Z\times{\mathfrak{S}}_{m-N}$, where $Z$ is the centralizer of
$\{\phi(T_{i})\mid i=1,3,4,5,6,0\}$ in ${\mathfrak{S}}_{N}$. By (3), the
support of each nontrivial cycle of $\varphi_{j}(T_{1})$ intersects
$S(\{\varphi_{j}(T_{i})\mid i=0,2,3,4,5,6\})$ nontrivially, and, by ($\ast$),
\begin{align*}
&  S(\{\varphi_{j}(T_{i})\mid i=0,2,3,4,5,6\})\\
=  &  S(\{\varphi_{1}(T_{i})\mid i=0,2,3,4,5,6\})\\
\subseteq &  \{1,\dots,N\}\,\text{.}%
\end{align*}
Thus, $\varphi_{j}(T_{1})\in Z$. Now, since $T_{j}^{\prime}$ and
$T_{1}^{\prime}$ are conjugate in ${\mathcal{M}}_{3,1}$, $\varphi_{j}(T_{1})$
and $\varphi_{1}(T_{1})$ share the same cycle decomposition type; hence, by
(2), $\varphi_{j}(T_{1})=\varphi_{1}(T_{1})$. To complete the proof, observe
that transitivity forces $m=N$.\qed

\begin{proposition}
\label{11107P4.3}Let $n\geq0$. Then $\mathcal{M}_{3,n}$ has precisely two
proper subgroups of index at most $36$ up to conjugation, namely,
$\mathcal{O}_{3,n}^{-}$ of index $28$, and $\mathcal{O}_{3,n}^{+}$ of index
$36$.
\end{proposition}

\noindent\textbf{Proof.} For $n\ge1$, the claim is a direct consequence of
Lemma \ref{11107L4.1} and the definitions of $\mathcal{O}_{3,n}^{-}$
respectively $\mathcal{O}_{3,n}^{+}$. By Theorem~\ref{11107T2.1},
${\mathcal{M}}_{3,0}$ is a quotient of ${\mathcal{M}}_{3,1}$ by one additional
relation, so for the case $n=0$ it is sufficient to check that the
representations of ${\mathcal{M}}_{3,1}$ given in Lemma \ref{11107L4.1}
satisfy the additional relation of ${\mathcal{M}}_{3,0}$; this is the case, as
can easily be verified. \qed


\pagestyle{myheadings} \markright{{Finite index subgroups of mapping class groups}  \quad\sc{\today}}

\part{Induction arguments\label{PartInduction}}

\pagestyle{myheadings} \markright{{Finite index subgroups of mapping class groups}  \quad\sc{\today}}

We turn now to the proof of our main result, Theorem \ref{11107T0.2}. As
pointed out before, we argue by induction on the genus. Recall that the case
$g=3$ is proved in Section \ref{11107S4} (see Proposition~\ref{11107P4.3}). Thus:

\begin{itemize}
\item \emph{from now on, we suppose that} $g\geq4$ \emph{plus the inductive
hypothesis} \emph{that Theorem \ref{11107T0.2} holds for a surface of genus}
$g-1$\emph{.}
\end{itemize}

Recall that we have defined $N_{g}^{-}=2^{g-1}(2^{g}-1)$ and $N_{g}%
^{+}=2^{g-1}(2^{g}+1)$. Throughout the arguments below, we shall rely on the
following numerical relationships for $g\geq4$.
\begin{gather*}
N_{g}^{+}=3\,N_{g-1}^{+}+N_{g-1}^{-}\,,\quad N_{g}^{-}=3\,N_{g-1}^{-}%
+N_{g-1}^{+}\,,\\
4\,N_{g-1}^{-}<N_{g}^{-}<N_{g}^{+}<5\,N_{g-1}^{-}<2\,N_{g}^{-}\,.
\end{gather*}
(Actually, only the third inequality requires $g\geq4$; the others hold for
$g\geq2$.)

Theorem \ref{11107T0.1} will be entirely proved in Section \ref{11107S7}.
Theorem \ref{11107T0.3} will be proved in Sections \ref{11107S8}
and~\ref{11107S9}.

\section{Factorization through symplectic groups\label{11107S5}}

Our goal in this section is to prove the following theorem.

\begin{theorem}
\label{11107T5.1}For $g\geq4$, let $\varphi:{\mathcal{M}}_{g,1}\rightarrow
{\mathfrak{S}}_{m}$ be a nontrivial homomorphism, with $m<5\,N_{g-1}^{-}$.
Then there is a homomorphism $\bar{\varphi}:\mathrm{Sp}_{2g}({\mathbb{F}}%
_{2})\rightarrow{\mathfrak{S}}_{m}$ such that the following diagram commutes.
\[
\xymatrix{
\MM_{g,1} \ar[d]_{\theta_{g,1}} \ar[dr]^{\varphi} \\
\mathrm{Sp}_{2g}(\F_2) \ar[r]_{\bar\varphi} & \SSS_m}
\]

\end{theorem}

\bigskip\noindent The proof strategy is as follows. To show that the kernel of
$\varphi$ contains the kernel of $\theta_{g,1}$, we first prove that the image
$w_{i}=\varphi(T_{i})$ has order $2$ in ${\mathfrak{S}}_{m}$. For this, we
consider the cycle decomposition type $(1)^{\ell_{1}}(2)^{\ell_{2}}%
\cdots(m)^{\ell_{m}}$ of the permutation $w_{i}$. We first exclude the
possibility of cycles of length at least $5$, then of cycles of length $4$,
and finally -- the most delicate case -- we exclude cycles of length $3$.
Hence, $w_{i}$ is reduced to being an involution after all. It then remains to
show that the kernel of $\theta_{g,1}$ is the normal closure of the square of
any of our standard generators of~${\mathcal{M}}_{g,1}$.

\subsection{Reduction to involutions\label{11107S5.1}}

\begin{proposition}
\label{11107P5.2}Let $g\geq4$, and let $b$ be a nonseparating simple closed
curve in $\Sigma_{g,1}$. Then the centralizer $\mathcal{Z}_{b}$ of $T_{b}$ in
$\mathcal{M}_{g,1}$ contains an index $2$ subgroup $\mathcal{Z}_{b}^{+}$ with
the properties:

\begin{enumerate}
\item[(a)] $\mathcal{Z}_{b}^{+}$ is perfect;

\item[(b)] $N_{g-1}^{-}\leq\mathrm{mi}(\mathcal{Z}_{b}^{+})$; and

\item[(c)] $T_{b}\in\mathcal{Z}_{b}^{+}$.
\end{enumerate}
\end{proposition}

\noindent\textbf{Proof.} It is known that ${\mathcal{Z}}_{b}$ is the set of
mapping classes that fix the curve $b$ up to isotopy (see \textsl{e.g.}
\cite{ParRol1}). We take ${\mathcal{Z}}_{b}^{+}$ to be the subgroup of
${\mathcal{Z}}_{b}$ consisting of those classes that also preserve the
orientation of $b$. Evidently, ${\mathcal{Z}}_{b}^{+}$ has index $2$ in
${\mathcal{Z}}_{b}$ and contains $T_{b}$. From \cite{ParRol1}, ${\mathcal{Z}%
}_{b}^{+}$ is the image of ${\mathcal{M}}_{g-1,3}$ in ${\mathcal{M}}_{g,1}$
under the homomorphism induced by a quotient map $\Sigma_{g-1,3}%
\rightarrow\Sigma_{g,1}$ identifying two boundary circles with the closed
curve $b$. Since by \cite{Harer1} ${\mathcal{M}}_{g-1,3}$ is perfect, we have
(a). From the epimorphism ${\mathcal{M}}_{g-1,3}\rightarrow{\mathcal{Z}}%
_{b}^{+}$ we know that $\mathrm{mi}({\mathcal{M}}_{g-1,3})\leq\mathrm{mi}%
({\mathcal{Z}}_{b}^{+})$, and by induction we have $\mathrm{mi}({\mathcal{M}%
}_{g-1,3})=N_{g-1}^{-}$, thus $N_{g-1}^{-}\leq\mathrm{mi}({\mathcal{Z}}%
_{b}^{+})$. \qed

\bigskip\noindent As ever, we consider the simple closed curves $a_{0}%
,a_{1},\dots,a_{2g+1}$ illustrated in Figure 2.1, we denote by $T_{i}$ the
Dehn twist about $a_{i}$, and write $w_{i}$ for the image under $\varphi$ of
the Dehn twist $T_{i}$ ($i=0,1,\dots,2g+1$). In the sequel we repeatedly use
the fact that, since $T_{i}$ and $T_{j}$ are conjugate in the mapping class
group, $w_{i}=\varphi(T_{i})$ and $w_{j}=\varphi(T_{j})$ are conjugate in the
symmetric group, and so share the same cycle decomposition type, say
$(1)^{\ell_{1}}(2)^{\ell_{2}}\cdots(m)^{\ell_{m}}$, where $\sum_{k=1}^{m}
k\ell_{k}=m$.

\noindent The fact that $T_{i}\in{\mathcal{Z}}_{a_{i}}^{+}$ implies that,
whenever $\ell_{k}>0$ with $k>1$, $\varphi({\mathcal{Z}}_{a_{i}}^{+})$ acts
nontrivially on the union of the $k$-orbits of $w_{i}$. Therefore, the above
proposition combines with Lemma~\ref{11107L1.1}\thinspace(3)\thinspace
(a),\thinspace(c) to yield the following.

\begin{corollary}
\label{11107C5.3}For $g\geq4$, let $\varphi:\mathcal{M}_{g,1}\rightarrow
\mathfrak{S}_{m}$ be a nontrivial homomorphism, with $m<5N_{g-1}^{-}$. Then
there exists $k\in\{2,3,4\}$ such that every $w_{i}$ has the same cycle
decomposition type $(1)^{\ell_{1}}(k)^{\ell_{k}}$ with $\ell_{k}\geq
N_{g-1}^{-}$. \qed

\end{corollary}

\begin{lemma}
\label{11107L5.4}For $g\geq3$, let $\varphi:\mathcal{M}_{g,1}\rightarrow
\mathfrak{S}_{m}$ be a group homomorphism. If $S(w_{2g-2})\cap S(w_{2g}%
)=\emptyset$, then $\varphi$ is trivial.
\end{lemma}

\noindent\textbf{Proof.} First note that, since conjugate permutations have
bijective supports, by Lemma \ref{11107L2.2} the cardinality $r=|S(w_{i})|$ is
independent of choice of $i$, and for all $3$-chains $(a_{i},a_{j},a_{k})$
yields the same cardinalities $s=|S(w_{i})\cap S(w_{j})|$ and $t=|S(w_{i})\cap
S(w_{k})|$. By assumption, the $3$-chain $(a_{2g-2},a_{2g-1},a_{2g})$ gives
$t=0$, which implies that
\[
(S(w_{2g-2})\cap S(w_{2g-1}))\sqcup(S(w_{2g-1})\cap S(w_{2g}))\subseteq
S(w_{2g-1})\,,
\]
and so $r\geq2s$. On the other hand, Lemma \ref{11107L1.2} asserts that
$r\leq2s$. Thus, for any $3$-chain $(a_{i},a_{j},a_{k})$ we have
\[
S(w_{j})=(S(w_{i})\cap S(w_{j}))\sqcup(S(w_{j})\cap S(w_{k}))\,.
\]
Turning now to the $3$-chains $(a_{0},a_{4},a_{3})$ and $(a_{0},a_{4},a_{5})$
and $(a_{3},a_{4},a_{5})$, we observe that, since $t=0$, the subsets
$S(w_{0})\cap S(w_{4})$, $S(w_{3})\cap S(w_{4})$ and $S(w_{4})\cap S(w_{5})$
of $S(w_{4})$ are pairwise disjoint, but each of cardinality half that of
$S(w_{4})$. For avoidance of contradiction, it must be that $r=s=0$; whence,
since by Theorem \ref{11107T2.1} the image of $\varphi$ is generated by the
$w_{i}$, $\varphi$ is trivial. \qed

\begin{proposition}
\label{11107P5.5}For $g\geq4$, let $\varphi:\mathcal{M}_{g,1}\rightarrow
\mathfrak{S}_{m}$ with $m<5N_{g-1}^{-}$. Then for $i=0,1,\ldots,2g+1$ there
are no $4$-cycles in the cycle decomposition of $w_{i}$.
\end{proposition}

\noindent\textbf{Proof.} Again write $\ell_{k}$ for the number of $k$-cycles
in the decomposition of $w_{2g}$, and suppose that $\ell_{4}>0$. Then
Corollary \ref{11107C5.3} forces the cycle decomposition type of every $w_{j}$
to be $(1)^{\ell_{1}}(4)^{\ell_{4}}$ where $\ell_{4}\geq N_{g-1}^{-}$. The
image $\varphi({\mathcal{M}}_{g,1})$ contains the subgroup $H=\langle
w_{2g},w_{2g+1}\rangle$, which, since $w_{2g}$ is a product of $4$-cycles, is
a homomorphic image (via $x_{i}\mapsto w_{2g+i}$) of the group
\[
\widehat{H}=\langle x_{0},x_{1}\,\mid\,x_{0}^{4}=x_{1}^{4}=1,\ x_{0}x_{1}%
x_{0}=x_{1}x_{0}x_{1}\rangle
\]
of order $96$ (denoted $\langle-2,3\,|\,4\rangle$ in \cite[p.74]{CoxMos1}).

\noindent From now on ${\mathcal{M}}_{g-1,1}$ is considered as the subgroup of
${\mathcal{M}}_{g,1}$ generated by $T_{0},T_{1},\dots,T_{2g-2}$. We write
$\widehat{\Omega}_{k}$ for the union of the orbits $\Omega_{k,i}$
($i=1,\dots,h_{k}$) of cardinality $k$ of $H$. Thus, $|\widehat{\Omega}%
_{k}|=kh_{k}$ and $k$ divides $96$. From the fact that the generators $w_{2g}$
and $w_{2g+1}$ of $H$ contain no nontrivial cycles of length less than $4$ it
follows that $h_{2}=h_{3}=0$. Since $\varphi({\mathcal{M}}_{g-1,1})$ lies in
the centralizer of $H$, it acts, for any $k$, both on $\widehat{\Omega}_{k}$
and on the set $\{ \Omega_{k,i} \mid i=1,\dots,h_{k} \}$.
The former corresponds to a homomorphism $\psi_{k}:{\mathcal{M}}%
_{g-1,1}\rightarrow{\mathfrak{S}}_{kh_{k}}$ and the latter to a homomorphism
$\nu_{k}:{\mathcal{M}}_{g-1,1}\rightarrow{\mathfrak{S}}_{h_{k}}$.

\medskip\noindent\textbf{Claim.} \emph{ If $\nu_{k}$ is trivial, then so is
$\psi_{k}$.}

\medskip\noindent To prove this, since $\nu_{1}=\psi_{1}$ we may assume that
$k\geq4$. Assuming that $\nu_{k}$ is trivial, we observe that each $k$-orbit
$\Omega_{k,i}$ of $H$ must be invariant under the action of $\varphi
({\mathcal{M}}_{g-1,1})$. Now partition $\Omega_{k,i}$ by means of the orbits
of $w_{2g}$, as
\[
\Omega_{k,i}={\mathcal{P}}_{1}\sqcup\cdots\sqcup{\mathcal{P}}_{p}\sqcup F\,,
\]
where each orbit ${\mathcal{P}}_{j}$ of $w_{2g}$ has length $4$, and $F$
consists of the fixed points of $w_{2g}$. Then $\varphi({\mathcal{M}}%
_{g-1,1})$ acts on ${\mathcal{P}}_{1}\sqcup\cdots\sqcup{\mathcal{P}}_{p}$, on
$\{{\mathcal{P}}_{1},\cdots,{\mathcal{P}}_{p}\}$, and on $F$. The action on
the set $\{{\mathcal{P}}_{1},\dots,{\mathcal{P}}_{p}\}$ must be trivial
because
\[
p\leq\frac{k}{4}\leq\frac{96}{4}<2^{g-2}(2^{g-1}-1)=N_{g-1}^{-}=\mathrm{mi}%
({\mathcal{M}}_{g-1,1})
\]
(by induction). This leaves $\varphi({\mathcal{M}}_{g-1,1})$ acting on each
$4$-orbit ${\mathcal{P}}_{j}$, where the action must again be trivial, since
certainly $4<\mathrm{mi}({\mathcal{M}}_{g-1,1})$. Hence, the action of
$\varphi({\mathcal{M}}_{g-1,1})$ on $\Omega_{k,i}$ stabilises $\Omega_{k,i}
\cap S(w_{2g})$ pointwise, and likewise $\Omega_{k,i} \cap S(w_{2g+1})$.
However, since $H$ is generated by $w_{2g}$ and $w_{2g+1}$, the union of
$\Omega_{k,i} \cap S(w_{2g})$ and $\Omega_{k,i} \cap S(w_{2g+1})$ is equal to
$\Omega_{k,i}$. Thus, the action of $\varphi({\mathcal{M}}_{g-1,1})$ is
trivial on each $\Omega_{k,i}$, and therefore also on their union
$\widehat{\Omega}_{k}$ as required.

\smallskip\noindent In consequence, if $\varphi({\mathcal{M}}_{g-1,1})$ acts
nontrivially on some $\widehat{\Omega}_{k}$, then $h_{k}\geq\mathrm{mi}%
({\mathcal{M}}_{g-1,1})$. Since
\[
kh_{k}=|\widehat{\Omega}_{k}|\leq m<5\,N_{g-1}^{-}=5\,\mathrm{mi}%
({\mathcal{M}}_{g-1,1})\,,
\]
the only possibilities are $k=1$ or $k=4$. From Lemma \ref{11107L5.4},
${\mathcal{M}}_{g-1,1}$ must act nontrivially on some $\widehat{\Omega}_{k}$
with $k>1$, and therefore $k=4$ indeed occurs.

\noindent The inequality $h_{4}\geq N_{g-1}^{-}$ implies that $h_{1}\leq
m-4h_{4}<N_{g-1}^{-}$, thus ${\mathcal{M}}_{g-1,1}$ acts trivially on
$\widehat{\Omega}_{1}$. It follows that
\[
S(\varphi({\mathcal{M}}_{g-1,1}))\subseteq\widehat{\Omega}_{4}\subseteq
S(H)\,,
\]
where, for a subgroup $G$ of ${\mathfrak{S}}_{m}$, $S(G)$ denotes its support.
Since $S(\varphi({\mathcal{M}}_{g-1,1}))$ contains $S(\langle w_{1}%
,w_{2}\rangle)$, which by Lemma \ref{11107L2.2} has the same cardinality as
$S(H)$, the above inclusions are indeed equalities.

\noindent So, $w_{2g}$ is a product of $4$-cycles and all nontrivial orbits of
$H=\langle w_{2g},w_{2g+1}\rangle$ have length~$4$. Applying again Lemma
\ref{11107L2.2}, we get that $w_{4}$ is a product of $4$-cycles and, for
$i\in\{0,3,5\}$, all nontrivial orbits of $\langle w_{4},w_{i}\rangle$ have
length $4$. By Lemma \ref{11107L1.5} we deduce that $w_{0}=w_{3}=w_{5}$. But,
since $w_{3}=w_{5}$, we also have
\[
w_{3}=(w_{5}w_{6})w_{3}(w_{5}w_{6})^{-1}=(w_{5}w_{6})w_{5}(w_{5}w_{6}%
)^{-1}=w_{6}\,\text{.}%
\]
Thus $w_{5}=w_{6}$, and therefore $w_{0}=w_{1}=\cdots=w_{2g}$. It follows that
the image of $\varphi$ is cyclic, and hence, because ${\mathcal{M}}_{g,1}$ is
perfect, $\varphi$ is trivial -- a contradiction. \qed

\begin{proposition}
\label{11107P5.6}For $g\geq4$, let $\varphi:{\mathcal{M}}_{g,1}\rightarrow
{\mathfrak{S}}_{m}$ with $m<5N_{g-1}^{-}$. Then for $0\leq i\leq2g+1$ there
are no $3$-cycles in the cycle decomposition of $w_{i}$.
\end{proposition}

\noindent\textbf{Proof.} Suppose to the contrary that $\ell_{3}>0$. In
particular, $\varphi$ cannot be trivial. Then Corollary~\ref{11107C5.3} forces
the cycle decomposition type of each $w_{i}$ to be $(1)^{\ell_{1}}%
(3)^{\ell_{3}}$ where $\ell_{3}\geq N_{g-1}^{-}$. The image $\varphi
({\mathcal{M}}_{g,1})$ contains the subgroup $H=\langle w_{2g},w_{2g+1}%
\rangle$, which, since $w_{2g}$ is a product of $3$-cycles, is a homomorphic
image (via $x_{i}\mapsto w_{2g+i}$) of the group
\[
\widehat{H}=\langle x_{0},x_{1}\,\mid\,x_{0}^{3}=x_{1}^{3}=1,\ x_{0}x_{1}%
x_{0}=x_{1}x_{0}x_{1}\rangle
\]
of order $24$. Again, we consider ${\mathcal{M}}_{g-1,1}$ as the subgroup of
${\mathcal{M}}_{g,1}$ generated by $T_{0},T_{1},\dots,T_{2g-2}$. Moreover, we
write $\widehat{\Omega}_{k}$ for the union of the orbits $\Omega_{k,i}$
($i=1,\dots,h_{k}$) of cardinality $k$ of $H$. Thus, $|\widehat{\Omega}%
_{k}|=kh_{k}$ and $k$ divides $24$. From the fact that the generators $w_{2g}$
and $w_{2g+1}$ of $H$ contain no nontrivial cycles of length less than $3$ it
follows that $h_{2}=0$.

\noindent Since $\varphi({\mathcal{M}}_{g-1,1})$ lies in the centralizer of
$H$, it acts, for any $k$, both on $\widehat{\Omega}_{k}$ and on the set $\{
\Omega_{k,i} \mid i=1,\dots,h_{k} \}$. The former corresponds to a
homomorphism $\psi_{k}:{\mathcal{M}}_{g-1,1}\rightarrow{\mathfrak{S}}_{kh_{k}%
}$ and the latter to a homomorphism $\nu_{k}:{\mathcal{M}}_{g-1,1}%
\rightarrow{\mathfrak{S}}_{h_{k}}$.

\medskip\noindent\textbf{Claim 1.} \emph{If $\nu_{k}$ is trivial, then so is
$\psi_{k}$.}

\medskip Indeed, suppose that $\nu_{k}$ is trivial. Then each $k$-orbit
$\Omega_{k,i}$ of $H$ must be invariant under the action of $\varphi
({\mathcal{M}}_{g-1,1})$. But
\[
|\Omega_{k,i}|=k\leq24<2^{g-2}(2^{g-1}-1)=\mathrm{mi}({\mathcal{M}}%
_{g-1,1})\,,
\]
thus $\varphi({\mathcal{M}}_{g-1,1})$ acts trivially on each $\Omega_{k,i}$,
that is, ${\mathcal{M}}_{g-1,1}$ acts trivially on $\widehat{\Omega}_{k}$.

\medskip\noindent\textbf{Claim 2.} \emph{$\psi_{k}$ is trivial if $k\geq5$.
Moreover, either $\psi_{3}$ or $\psi_{4}$ is nontrivial, but not both.}

\medskip Indeed, if $\psi_{k}$ is nontrivial, then $\nu_{k}$ is nontrivial by
Claim 1, thus
\[
\frac{m}{k}\geq h_{k}\geq\mathrm{mi}({\mathcal{M}}_{g-1,1})=N_{g-1}^{-}\,,
\]
therefore $k\leq4$ as $m<5\,N_{g-1}^{-}$. If both $\psi_{3}$ and $\psi_{4}$
are nontrivial, then
\[
m\geq3h_{3}+4h_{4}\geq7N_{g-1}^{-}\,,
\]
which also contradicts $m<5\,N_{g-1}^{-}$. Finally, by Lemma \ref{11107L5.4}
one of the $\psi_{k}$ with $k\geq3$ must be nontrivial, so either $\psi_{3}$
or $\psi_{4}$ is nontrivial.

\medskip\noindent\textbf{Claim 3.} $\psi_{3}$\emph{ is nontrivial.}

\medskip\noindent Suppose otherwise that $\psi_{3}$ is trivial. Then, by Claim
2, $\psi_{4}$ is nontrivial, thus $h_{4}\geq N_{g-1}^{-}$. As $h_{1}%
+4h_{4}\leq m<5N_{g-1}^{-}$, it follows that $h_{1}<N_{g-1}^{-}$, thus
${\mathcal{M}}_{g-1,1}$ acts trivially on $\widehat{\Omega}_{1}$. So,
\[
S(\varphi({\mathcal{M}}_{g-1,1}))\subseteq\widehat{\Omega}_{4}\subseteq
S(H)\,.
\]
Since $S(\varphi({\mathcal{M}}_{g-1,1}))$ contains $S(\langle w_{1}%
,w_{2}\rangle)$, which has the same cardinality as $S(H)$, the above
inclusions are indeed equalities.

\medskip\noindent So, all nontrivial orbits of $H$ have length $4$. Applying
Lemma \ref{11107L2.2}, we get that $w_{2}$ is a product of disjoint $3$-cycles
and, for $i\in\{1,3\}$, all nontrivial orbits of $\langle w_{2},w_{i}\rangle$
have length $4$. By Lemma~\ref{11107L1.4} it follows that $w_{1}=w_{3}$.
Hence, we also have
\[
w_{1}=(w_{3}w_{4})w_{1}(w_{3}w_{4})^{-1}=(w_{3}w_{4})w_{3}(w_{3}w_{4}%
)^{-1}=w_{4}\,,
\]
thus $w_{3}=w_{4}$, therefore $w_{0}=w_{1}=\cdots=w_{2g}$. This implies that
the image of $\varphi$ is cyclic, and hence, because ${\mathcal{M}}_{g,1}$ is
perfect, $\varphi$ is trivial -- a contradiction.

\medskip\noindent\textbf{Claim 4.} $\psi_{1}$\textit{ }\emph{is nontrivial.}

\medskip\noindent Suppose instead that $\psi_{1}$ is trivial. Then we have the
inclusions
\[
S(\varphi({\mathcal{M}}_{g-1,1}))\subseteq\widehat{\Omega}_{3}\subseteq S(H)
\]
which imply as in the previous claim that $S(H)=\widehat{\Omega}_{3}$. So, all
nontrivial orbits of $w_{2g}$ and $\langle w_{2g},w_{2g+1}\rangle$ have length
$3$; thus, by Lemma~\ref{11107L1.3}, $w_{2g}=w_{2g+1}$. It follows that the
image of $\varphi$ is cyclic, and hence, because ${\mathcal{M}}_{g,1}$ is
perfect, $\varphi$ is trivial -- a contradiction.

\medskip\noindent The following claim is a direct consequence of the braid
relation between $w_{2g}$ and $w_{2g+1}$ (see the proof of Lemma
\ref{11107L1.3}).

\medskip\noindent\textbf{Claim 5.} $\Omega_{3,1},\dots,\Omega_{3,h_{3}}%
$\textit{ }\emph{are precisely the common}\textit{ }$3$\emph{-orbits of
}$w_{2g}$\emph{ and}\textit{ }$w_{2g+1}$\emph{.}

\medskip We turn now to conclude the proof of Proposition \ref{11107P5.6} with
counting arguments. If $(b_{1},\dots,b_{p})$ is a $p$-chain of nonseparating
closed curves, we denote by $c_{p}(b_{1},\dots,b_{p})$ the number of common
$3$-orbits of $\varphi(T_{b_{1}}),\dots,\varphi(T_{b_{p}})$. We do not assume
that the whole chain is nonseparating, but we suppose that each pair
$(b_{i},b_{i+1})$ is nonseparating, so that $(T_{b_{i}},T_{b_{i+1}})$ is
conjugate to $(T_{2g},T_{2g+1})$ (see Lemma \ref{11107L2.2}). By Claims 1, 3
and 5 we have
\[
c_{2}(b_{1},b_{2})=c_{2}(a_{2g},a_{2g+1})=h_{3}\geq N_{g-1}^{-}\,,
\]
and by Claim 4 we have $h_{1}\geq N_{g-1}^{-}$. Recall also that $\ell_{3}$ is
the number of $3$-orbits of $w_{2g}$. Note that the supports of the $3$-orbits
of $w_{2g}$ that are not included in $\widehat{\Omega}_{3}$ are included in
$X=\{1,\dots,m\}-(\widehat{\Omega}_{3}\sqcup\widehat{\Omega}_{1})$; thus there
are at most $\left\vert X\right\vert /3$ of them. Again by Lemma
\ref{11107L2.2} it follows that for all $i\in\{2,\dots,p\}$ among the
$3$-orbits of $\varphi(T_{i})$ there are at most $\left\vert X\right\vert /3$
that are not $3$-orbits of $\varphi(T_{i-1})$. We therefore obtain in turn
\begin{gather*}
c_{3}(b_{1},b_{2},b_{3})\geq h_{3}-\frac{|X|}{3}\,,\quad c_{4}(b_{1}%
,b_{2},b_{3},b_{4})\geq h_{3}-\frac{2|X|}{3}\,,\\
c_{5}(b_{1},b_{2},b_{3},b_{4},b_{5})\geq h_{3}-\frac{3|X|}{3}=h_{3}-|X|\,.
\end{gather*}
Using again the assumption that $m<5N_{g-1}^{-}$, we have
\[
|X|=m-h_{1}-3h_{3}<4N_{g-1}^{-}-3h_{3}\leq h_{3}\,.
\]
Hence
\[
c_{5}(b_{1},b_{2},b_{3},b_{4},b_{5})\geq h_{3}-|X|>0\,.
\]

\medskip So, there is a common $3$-cycle in the decomposition of every
generator $w_{0},\dots,w_{4}$ of $\varphi({\mathcal{M}}_{2,1})$. By
restricting attention to the support of this $3$-cycle, we deduce that
$\varphi$ induces a nontrivial homomorphism from ${\mathcal{M}}_{2,1}$ to
${\mathfrak{S}}_{3}$ whose image is generated by elements of order $3$; in
other words, from ${\mathcal{M}}_{2,1}$ onto the cyclic group of order $3$.
This, however, contradicts Proposition~\ref{11107P3.3}.\qed


\subsection{Normal generation of $\mathrm{Ker}\,\theta_{g,1}$ by the square of
a Dehn twist\label{11107S5.2}}

Thanks to Corollary \ref{11107C5.3} and Propositions \ref{11107P5.5} and
\ref{11107P5.6}, we know (under the inductive hypothesis that Theorem
\ref{11107T0.2} holds for a surface of genus $g-1$) that, if $g\ge4$ and
$\varphi:{\mathcal{M}}_{g,1}\rightarrow{\mathfrak{S}}_{m}$ is a permutation
representation with $m<5\,N_{g-1}^{-}$, then $w_{i}=\varphi(T_{i})$ is an
involution for all $i\in\{0,1,\dots,2g+1\}$. So, in order to prove Theorem
\ref{11107T5.1}, still under the inductive hypothesis, it suffices to show the following.

\begin{theorem}
\label{11107T6.1}Let $g\geq2$. Then the kernel of $\theta_{g,1}$ is the
smallest normal subgroup of $\mathcal{M}_{g,1}$ containing $T_{1}^{2}$.
\end{theorem}

\medskip This theorem is known to experts but, as far as we know, is nowhere
in the literature. So, we include a proof but skip some details.

\medskip\noindent\textbf{Proof.} We say that a normal subgroup $H$ of a group
$G$ is \emph{normally generated} by a subset $S\subseteq H$ if $H$ is the
smallest normal subgroup of $G$ that contains $S$. For $i,j\in\{1,\dots,g\}$
we denote by $e_{i,j}$ the $g\times g$ matrix whose entries are all zero
except the $(i,j)$-th entry which is equal to $1$. Then the kernel of $\mu
_{g}:\mathrm{Sp}_{2g}({\mathbb{Z}})\rightarrow\mathrm{Sp}_{2g}({\mathbb{F}%
}_{2})$ is normally generated by $\left(
\begin{matrix}
I_{g} & 2\,e_{1,1}\\
0 & I_{g}%
\end{matrix}
\right)  $, where $I_{g}$ denotes the identity matrix. This fact can be
deduced from \cite{BaMiSe1} (see also \cite{Tits1}) and its proof is left to
the reader.

\noindent We denote by $\widehat{\theta}_{g}:{\mathcal{M}}_{g,1}%
\rightarrow\mathrm{Sp}_{2g}({\mathbb{Z}})$ the natural epimorphism, so that
$\theta_{g,1}=\mu_{g}\circ\widehat{\theta}_{g}$. Consider the simple closed
curves $c_{2},\dots,c_{g}$, $c_{2}^{\prime},\dots,c_{g}^{\prime}$,
$d_{2},\dots,d_{g}$ illustrated in Figure 5.1. By \cite{Powel1} and
\cite{Birma1} (see also \cite{Putma1}), $\mathrm{Ker}\,\widehat{\theta}_{g}$
is normally generated by $\{T_{c_{i}}T_{c_{i}^{\prime}}^{-1},T_{d_{i}}
\mid2\leq i\leq g\}$. On the other hand, a direct calculation shows that
$\widehat{\theta}_{g}(T_{1}^{2})=\left(
\begin{matrix}
I_{g} & 2e_{1,1}\\
0 & I_{g}%
\end{matrix}
\right)  $. By the above, it follows that $\mathrm{Ker}\,\theta_{g,1}$ is
normally generated by $\{T_{c_{i}}T_{c_{i}^{\prime}}^{-1},T_{d_{i}} \mid2\leq
i\leq g\}\cup\{T_{1}^{2}\}$.

\begin{figure}[tbh]
\bigskip\centerline{
\setlength{\unitlength}{0.5cm}
\begin{picture}(20.5,4)
\put(0,0){\includegraphics[width=10.25cm]{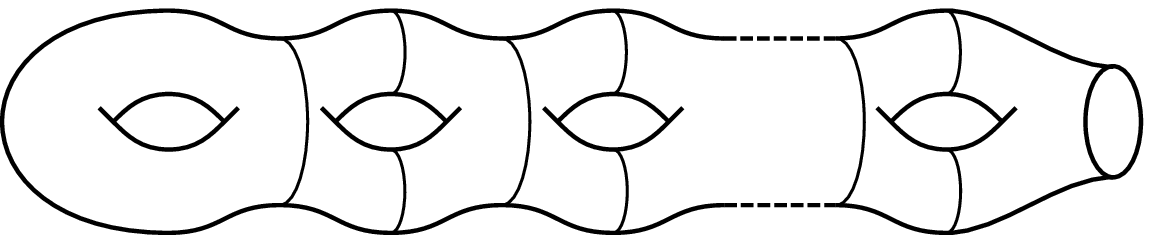}}
\put(4.9,3.9){\small $d_2$}
\put(8.9,3.9){\small $d_3$}
\put(14.9,3.9){\small $d_g$}
\put(7,-0.5){\small $c_2$}
\put(11,-0.5){\small $c_3$}
\put(17,-0.5){\small $c_g$}
\put(6.9,4.4){\small $c_2'$}
\put(10.9,4.4){\small $c_3'$}
\put(16.9,4.4){\small $c_g'$}
\end{picture}} \bigskip
\centerline{{\bf Figure 5.1.} Some curves in $\Sigma_{g,1}$.} \bigskip
\end{figure}

\medskip Let $H$ be the normal subgroup of ${\mathcal{M}}_{g,1}$ normally
generated by $T_{1}^{2}$, and let $\pi:{\mathcal{M}}_{g,1}\rightarrow
{\mathcal{M}}_{g,1}/H$ be the quotient map. The following relations hold in
${\mathcal{M}}_{g,1}/H$:
\[%
\begin{array}
[c]{cl}%
\pi(T_{i})^{2}=1 & \quad\text{if }1\leq i\leq2g\,,\\
\noalign{\smallskip}\pi(T_{i})\pi(T_{j})=\pi(T_{j})\pi(T_{i}) & \quad\text{if
}|i-j|\geq2\,,\\
\noalign{\smallskip}\pi(T_{i})\pi(T_{i+1})\pi(T_{i})=\pi(T_{i+1})\pi(T_{i}%
)\pi(T_{i+1}) & \quad\text{if }1\leq i\leq2g-1\,.
\end{array}
\]
Thus, there exists a homomorphism $q:{\mathfrak{S}}_{2g+1}\rightarrow
{\mathcal{M}}_{g,1}/H$ that sends $s_{i}=(i,i+1)$ to $\pi(T_{i})$ for all
$1\leq i\leq2g$. On the other hand, the following relations hold in
${\mathcal{M}}_{g,1}$ (see \cite{LabPar1}, for instance):
\begin{gather*}
T_{d_{i}}=(T_{1}T_{2}\cdots T_{2i-2})^{4i-2}\,,\\
T_{c_{i}}T_{c_{i}^{\prime}}=(T_{1}T_{2}\cdots T_{2i-1})^{2i}\,,
\end{gather*}
whenever $i\in\{2,\dots,g\}$. Moreover, since $s_{1}s_{2}\cdots s_{j}$ is a
cycle of length $j+1$,
\[
(s_{1}s_{2}\cdots s_{j})^{j+1}=1\,,\quad\text{for }1\leq j\leq2g\,.
\]
It follows that
\begin{gather*}
\pi(T_{d_{i}})=q((s_{1}s_{2}\cdots s_{2i-2})^{4i-2})=q(1)=1\,,\\
\pi(T_{c_{i}}T_{c_{i}^{\prime}}^{-1})=\pi(T_{c_{i}}T_{c_{i}^{\prime}%
})=q((s_{1}s_{2}\cdots s_{2i-1})^{2i})=q(1)=1\,,
\end{gather*}
for all $i\in\{2,\dots,g\}$. Thus, $H=\mathrm{Ker}\,\theta_{g,1}$. \qed

\bigskip\noindent\textbf{Proof of Theorem \ref{11107T5.1}.} For $g\geq4$, let
$\varphi:{\mathcal{M}}_{g,1}\rightarrow{\mathfrak{S}}_{m}$ be a nontrivial
homomorphism, with $m<5\,N_{g-1}^{-}$. As pointed out before, thanks to
Corollary \ref{11107C5.3} and Propositions \ref{11107P5.5} and \ref{11107P5.6}%
, $w_{i}=\varphi(T_{i})$ is an involution for all $i\in\{0,1,\dots,2g+1\}$. By
Theorem \ref{11107T6.1} we conclude that there exists a homomorphism
$\bar{\varphi}:\mathrm{Sp}_{2g}({\mathbb{F}}_{2})\rightarrow{\mathfrak{S}}%
_{m}$ such that $\varphi=\bar{\varphi}\circ\theta_{g,1}$. \qed


\section{Large subgroups of $\mathrm{Sp}_{2g}({\mathbb{F}}_{2})$%
\label{11107S7}}

The aim of this section is to prove Theorem \ref{11107T0.1}, which, together
with Theorem \ref{11107T5.1}, proves Theorem \ref{11107T0.2} for the case of a
surface with $n=1$ boundary component. The extension to surfaces with several
boundary components will be the object of Section \ref{11107S9}.

\medskip The maximal subgroups of the finite classical groups have been
described using the classification of finite simple groups. (Below, $2^{n}$
denotes an elementary abelian group of rank $n$, $A\mathbf{\centerdot}B$ an
extension of a group $A$ by a group $B$, and $r$ a cyclic group of order $r$.)
Their orders are also known (see for example \cite{Taylo1}).

\begin{theorem}
\label{11107T7.1}\emph{(\cite[Theorem 4.1 and §3]{Liebe1}; \cite[§1]{Aschb1})}
If $H$ is a maximal subgroup of $\mathrm{Sp}_{2g}({\mathbb{F}}_{2})$ for
$g\geq1$, then one of the following holds.

\begin{itemize}
\item[\textbf{(1)}] $H=O_{2g}^{-}({\mathbb{F}}_{2})$, of order $2^{g^{2}%
-g+1}(2^{g}+1)\prod_{i=1}^{g-1}(2^{2i}-1)$.

\item[\textbf{(2)}] $H=O_{2g}^{+}({\mathbb{F}}_{2})$, of order $2^{g^{2}%
-g+1}(2^{g}-1)\prod_{i=1}^{g-1}(2^{2i}-1)$.

\item[\textbf{(3)}] $|H| \le2^{6g}$.

\item[\textbf{(4)}] $H \le{\mathfrak{S}}_{2g+2}$, of order at most $(2g+2)!$.

\item[\textbf{(5)}] $H=\mathrm{Sp}_{2k}({\mathbb{F}}_{2^{r}}%
)\mathbf{\centerdot}r$, where $r>1$ is a prime divisor of $g$ and $kr=g$.
Here, $H$ has order $r\cdot2^{\frac{g^{2}}{r}}\prod_{i=1}^{k}(2^{2ri}-1)$.

\item[\textbf{(6)}] $H=\mathrm{Sp}_{2k}({\mathbb{F}}_{2})\wr{\mathfrak{S}}%
_{r}$, where $r>1$ is a divisor of $g$ and $kr=g$. Here, $H$ has order
$r!\cdot2^{\frac{g^{2}}{r}}\prod_{i=1}^{k}(2^{2i}-1)^{r}$.

\item[\textbf{(7)}] $H$ is the stabilizer of a totally isotropic subspace of
${\mathbb{F}}_{2}^{2g}$ under the natural action of $\mathrm{Sp}%
_{2g}({\mathbb{F}}_{2})$. That is, $H$ is
\[
(2^{\frac{k(k+1)}{2}}\mathbf{\centerdot}2^{2k(g-k)})\rtimes(\mathrm{Sp}%
_{2g-2k}({\mathbb{F}}_{2})\times\mathrm{GL}_{k}({\mathbb{F}}_{2}))
\]
for some integer $k\in\{1,\dots,g\}$, and has order
\[
2^{g^{2}}\cdot\left(  \prod_{i=1}^{g-k}(2^{2i}-1)\right)  \cdot\left(
\prod_{i=1}^{k-1}(2^{i+1}-1)\right)  \,.
\]

\item[\textbf{(8)}] $H$ is the stabilizer of a nonsingular subspace of
${\mathbb{F}}_{2}^{2g}$ under the natural action of $\mathrm{Sp}%
_{2g}({\mathbb{F}}_{2})$. That is, $H$ is $\mathrm{Sp}_{2k}({\mathbb{F}}%
_{2})\times\mathrm{Sp}_{2g-2k}({\mathbb{F}}_{2})$ for some $k\in
\{1,\dots,g-1\}$, and has order
\[
2^{g^{2}+2k^{2}-2gk}\cdot\left(  \prod_{i=1}^{g-k}(2^{2i}-1)\right)
\cdot\left(  \prod_{i=1}^{k}(2^{2i}-1)\right)  \,.
\]

\end{itemize}

Moreover, the maximal subgroups falling within case (1) lie in a single
conjugacy class, as do the subgroups falling within case (2).
\end{theorem}

\bigskip\noindent\textbf{Proof.} The first part is shown in \cite[Theorem 4.1
and §3]{Liebe1}. Their orders may be found in \cite[pp.\thinspace19, 70,
141]{Taylo1}. For the last two cases, the description of $H$ follows from
\cite[p.\,63]{Wilso1}, for example. Cases 2 and 1 correspond to the Aschbacher
class ${\mathcal{C}}_{8}$ \cite[§1]{Aschb1}, that is, contain the subgroups of
$\mathrm{Sp}_{2g}({\mathbb{F}}_{2})$ stabilizing some quadratic form
polarizing to the symplectic bilinear form defining $\mathrm{Sp}%
_{2g}({\mathbb{F}}_{2})$; Case 2 corresponds to forms of Witt index $g$, Case
1 to forms of Witt index $g-1$. By \cite[Theorems~$B\Delta$ and $BO$]{Aschb1},
the action of $\mathrm{Sp}_{2g}({\mathbb{F}}_{2})$ is transitive on
stabilizers of forms of the same similarity type, showing that the groups
falling within Case 2, respectively the groups falling within Case 1, form a
single conjugacy class. \qed

\begin{corollary}
\label{11107C7.2}Let $g\geq3$ and let $H$ be a subgroup of $\mathrm{Sp}%
_{2g}(\mathbb{F}_{2})$ of index (strictly) less than $2\,N_{g}^{-}$. Then
either $H\cong O_{2g}^{-}(\mathbb{F}_{2})$ or $H\cong O_{2g}^{+}%
(\mathbb{F}_{2})$.
\end{corollary}

\bigskip\noindent\textbf{Remark.} $SO_{2g}^{-}({\mathbb{F}}_{2})$ is an index
$2$ subgroup of $O_{2g}^{-}({\mathbb{F}}_{2})$ and, therefore, an index $2\,
N_{g}^{-}$ subgroup of $\mathrm{Sp}_{2g}({\mathbb{F}}_{2})$.

\bigskip\noindent\textbf{Proof.} Recall again that $N_{g}^{-}=[\mathrm{Sp}%
_{2g}({\mathbb{F}}_{2}):O_{2g}^{-}({\mathbb{F}}_{2})]$ and note that the
result holds for proper subgroups of $O_{2g}^{-}({\mathbb{F}}_{2})$ and
$O_{2g}^{+}({\mathbb{F}}_{2})$ because $|O_{2g}^{+}({\mathbb{F}}_{2}%
)|<|O_{2g}^{-}({\mathbb{F}}_{2})|$. Thus, we need consider only maximal
subgroups $H$. For applications, the assertion is presented in terms of the
index; however, the data relate to the order, so one needs to check that of
the maximal subgroups listed in Theorem \ref{11107T7.1}, in all except for the
first two cases the order of $H$ is less than
\[
\frac{|O_{2g}^{-}({\mathbb{F}}_{2})|}{2}=2^{g^{2}-g}(2^{g}+1)\prod_{i=1}%
^{g-1}(2^{2i}-1)\,.
\]

\medskip For $g=3$ the result holds by \cite{CCNPW1}. (Alternatively, we can
compute the conjugacy classes of maximal subgroups of $\mathrm{Sp}%
_{6}({\mathbb{F}}_{2})$ using \textsc{Magma} \cite{BoCaPl1}: their orders are
1512, 4320, 4608, 10752, 12096, 23040, 40320, 51840. The maximal subgroups of
order 40320 and 51840 can be checked to be isomorphic to $O_{6}^{+}%
({\mathbb{F}}_{2})$ respectively $O_{6}^{-}({\mathbb{F}}_{2})$, so the claim holds.)

\medskip So, assume that $g\geq4$. In general, it is a routine matter to use
the formulae of Theorem \ref{11107T7.1} to check that $H$ has order less than
$|O_{2g}^{-}({\mathbb{F}}_{2})|/2$. To indicate the flavor of the
verification, we discuss the two most delicate cases, namely (6) and (8).

\bigskip\noindent Case (6). First observe that $r! < 2^{(r^{2}-2r+3)/2}$, and,
for $r$ in the range $[2,g]$, the function $(r^{2}-2r+3)/2+ \frac{g^{2}}{r}$
achieves its maximum value of $(g^{2}+3)/2$ at the endpoints. Thus, on putting
$j=ir$, we have that $|H| = r! \cdot2^{\frac{g^{2}}{r}} \prod_{i=1}^{k}
(2^{2i}-1)^{r}$ is bounded above by
\[
2^{(g^{2}+3)/2} \prod_{r|j,\ j\le g} (2^{2j}-1) <2^{(g^{2}+2g+3)/2} (2^{g}+1)
\prod_{r|j,\ j\le g-1} (2^{2j}-1)\,.
\]
Since $2^{(g^{2}+2g+3)/2} < 2^{g^{2}-g}(2^{2}-1)$, and $(2^{2}-1)
\prod_{r|j,\ j\le g-1} (2^{2j}-1) < \prod_{i=1}^{g-1}(2^{2i}-1)$, we obtain
$|H| < |O_{2g}^{-}({\mathbb{F}}_{2})| /2$, as required.

\bigskip\noindent Case (8). Since, for $k$ in the range $[1,g-1]$, the
function $2k(k-g)$ achieves its maximum of $-2(g-1)$ at the endpoints, we have
$2^{g^{2}+2k^{2}-2gk}<2^{g^{2}-g}$. Meanwhile, using the symmetry of the
product below to assume that $k\le\frac{g}{2}$, and observing that $i<i+g-k=j$
say, gives
\begin{align*}
&  \left(  \prod_{i=1}^{g-k}(2^{2i}-1)\right)  \cdot\left(  \prod_{i=1}%
^{k}(2^{2i}-1)\right) \\
=  &  (2^{2k}-1)\left(  \prod_{i=1}^{g-k}(2^{2i}-1)\right)  \cdot\left(
\prod_{i=1}^{k-1}(2^{2i}-1)\right) \\
<  &  (2^{g}+1)\left(  \prod_{i=1}^{g-k}(2^{2i}-1)\right)  \cdot\left(
\prod_{j=g-k+1}^{g-1}(2^{2j}-1)\right)  \,.
\end{align*}
When combined with the exponential inequality above, this again yields $|H| <
|O_{2g}^{-}({\mathbb{F}}_{2})|/2$. \qed

\bigskip\noindent\textbf{Proof of Theorem \ref{11107T0.1}.} As stated before,
we have $[\mathrm{Sp}_{2g}({\mathbb{F}}_{2}):O_{2g}^{-}({\mathbb{F}}%
_{2})]=N_{g}^{-}=2^{g-1}(2^{g}-1)$ and $[\mathrm{Sp}_{2g}({\mathbb{F}}%
_{2}):O_{2g}^{+}({\mathbb{F}}_{2})]=N_{g}^{+}=2^{g-1}(2^{g}+1)$ by
\cite[pp.\,70 and 141]{Taylo1}. Finally, Corollary~\ref{11107C7.2} ensures
that the conjugacy classes of $O_{2g}^{+}({\mathbb{F}}_{2})$ and $O_{2g}%
^{-}({\mathbb{F}}_{2})$ are the only conjugacy classes of subgroups of
$\mathrm{Sp}_{2g}({\mathbb{F}}_{2})$ of index less than $2\,N_{g}^{-}$. \qed

\bigskip\noindent Now, thanks to Theorem \ref{11107T5.1} and Theorem
\ref{11107T0.1}, we can prove our main theorem for the case of a surface with
$n=1$ boundary component. However, we have to keep the inductive hypothesis
(that Theorem \ref{11107T0.2} holds for a surface of genus $g-1$), as we have
used the equality $\mathrm{mi}({\mathcal{M}}_{g-1,3})=N_{g-1}^{-}$ in the
proof of Proposition \ref{11107P5.2}.

\begin{proposition}
\label{11107P7.3}Let $g\geq4$.

\begin{itemize}
\item[\textbf{(1)}] ${\mathcal{O}}_{g,1}^{-}$ is the unique subgroup of
${\mathcal{M}}_{g,1}$ of index $N_{g}^{-}=2^{g-1}(2^{g}-1)$, up to conjugation.

\item[\textbf{(2)}] ${\mathcal{O}}_{g,1}^{+}$ is the unique subgroup of
${\mathcal{M}}_{g,1}$ of index $N_{g}^{+}=2^{g-1}(2^{g}+1)$, up to conjugation.

\item[\textbf{(3)}] All the other proper subgroups of ${\mathcal{M}}_{g,1}$
are of index at least $5N_{g-1}^{-} > N_{g}^{+}$.
\end{itemize}
\end{proposition}

\medskip\noindent\textbf{Proof.} Let $H$ be a subgroup of ${\mathcal{M}}%
_{g,1}$ of index $m<5\,N_{g-1}^{-}$. By Theorem~\ref{11107T5.1} there is a
subgroup $\bar{H}$ of $\mathrm{Sp}_{2g}({\mathbb{F}}_{2})$ of index $m$ such
that $H=\theta_{g,1}^{-1}(\bar{H})$. Since $5\,N_{g-1}^{-}<2\,N_{g}^{-}$, by
Theorem \ref{11107T0.1}, we must necessarily have either $\bar{H}=O_{2g}%
^{-}({\mathbb{F}}_{2})$ or $\bar{H}=O_{2g}^{+}(\mathbb{F}_{2})$ up to
conjugation, thus either $H={\mathcal{O}}_{g,1}^{-}$ or $H={\mathcal{O}}%
_{g,1}^{+}$ up to conjugation. \qed


\section{Small symplectic representations of $\mathcal{M}_{g,n}$%
\label{11107S8}}

The aim of this section is to prove Proposition~\ref{11107P8.1}, which is the
same as Theorem~\ref{11107T0.3}, except that we do not state that the
decompositions in (1) are unique. The uniqueness will follow from Theorem
\ref{11107T0.2} proved in~Section \ref{11107S9}.

The decompositions of $\phi_{g,n}^{\pm}$ in (1) reflect analogous properties
of the corresponding representations of $\mathrm{Sp}_{2g}({\mathbb{F}}_{2})$,
which in turn arise in a geometric way; while this result may be more or less
known to experts, we could not locate a reference and hence establish it in
Lemma~\ref{L:SymplecticPresentations}. \medskip

We choose a basis $\{{\mathbf{e}}_{1},\dots,{\mathbf{e}}_{2g}\}$ of
$V={\mathbb{F}}_{2}^{2g}$ such that $({\mathbf{e}}_{i},{\mathbf{e}}_{2g+1-i})$
for $i=1,\dots,g$ are the symplectic pairs for the action of $\mathrm{Sp}%
_{2g}({\mathbb{F}}_{2})$; we write $\bar{i}$ for $2g+1-i$ to shorten notation.
For a vector ${\mathbf{x}}\in V$, we denote the components of ${\mathbf{x}}$
with respect to the basis $\{{\mathbf{e}}_{1},\dots,{\mathbf{e}}_{2g}\}$ by
$x_{1},\dots,x_{2g}$. We consider $\mathrm{Sp}_{2g-2}({\mathbb{F}}_{2})$ as a
subgroup of $\mathrm{Sp}_{2g}({\mathbb{F}}_{2})$ via the the standard
embedding defined by
\[
\omega\mapsto\left[
\begin{array}
[c]{c|ccc|c}%
1 & 0 & \cdots & 0 & 0\\\hline
\rule[5pt]{0pt}{4.5pt}0 &  &  &  & 0\\[-1ex]%
\vdots &  & \omega &  & \vdots\\
0 &  &  &  & 0\\\hline
\rule[5pt]{0pt}{4.5pt}0 & 0 & \cdots & 0 & 1
\end{array}
\right]
\]
with respect to the above basis. \medskip

For $g\ge2$ and $\epsilon\in\{\pm\}$, let $\overline{\phi}_{g}^{\,\epsilon}:
\mathrm{Sp}_{2g}({\mathbb{F}}_{2})\rightarrow{\mathfrak{S}}_{N_{g}^{\epsilon}%
}$ denote the permutation representation induced by the action of
$\mathrm{Sp}_{2g}({\mathbb{F}}_{2})$ on the right cosets of $O_{2g}^{\epsilon
}({\mathbb{F}}_{2})$ (cf.~Corollary~\ref{11107C7.2}).

\begin{lemma}
\label{L:SymplecticPresentations} Let $g\geq3$. Then,
\[
\overline{\phi}_{g}^{\,+}|_{\mathrm{Sp}_{2g-2}({\mathbb{F}}_{2})}%
\cong(\overline{\phi}_{g-1}^{\,+})^{3}\oplus\overline{\,\phi}_{g-1}^{-}%
\quad\text{and}\quad\overline{\phi}_{g}^{\,-}|_{\mathrm{Sp}_{2g-2}%
({\mathbb{F}}_{2})}\cong(\overline{\phi}_{g-1}^{\,-})^{3}\oplus\overline{\phi
}_{g-1}^{\,+}\,.
\]

\end{lemma}

\bigskip\noindent\textbf{Proof.} It is well-known (see, for example,
\cite{Dye1} or \cite[ch. 11]{Taylo1}) that $\overline{\phi}_{g}^{\,+}$ and
$\overline{\phi}_{g}^{\,-}$ arise from the action of $\mathrm{Sp}%
_{2g}({\mathbb{F}}_{2})$ on the two orbits of quadratic forms polarizing to
the symplectic form $\langle-,\,-\rangle$ preserved by $\mathrm{Sp}%
_{2g}({\mathbb{F}}_{2})$. The two orbits consist of the quadratic forms of
Witt index $g$ (type $+$) respectively Witt index $g-1$ (type $-$).

Let ${\mathcal{Q}}$ denote the set of quadratic forms polarizing to
$\langle-,\,-\rangle$. The elements of ${\mathcal{Q}}$ are of the form
\[
Q_{{\mathbf{b}}}({\mathbf{x}}) = \sum_{i=1}^{g} x_{i} x_{\bar{i}} +
\langle{\mathbf{x}},{\mathbf{b}}\rangle^{2} = \sum_{i=1}^{g} x_{i} x_{\bar{i}}
+ \langle{\mathbf{x}},{\mathbf{b}}\rangle\quad\mbox{ for $\bb\in V$.}
\]
The symplectic group $\mathrm{Sp}_{2g}({\mathbb{F}}_{2})$ acts on
${\mathcal{Q}}$ via $(\omega\cdot Q)({\mathbf{x}}) := Q(\omega{\mathbf{x}})$
for $\omega\in\mathrm{Sp}_{2g}({\mathbb{F}}_{2})$ and $Q\in{\mathcal{Q}}$. In
particular, $(\omega\cdot Q_{{\mathbf{b}}})({\mathbf{x}}) = Q_{{\mathbf{b}%
}^{\omega}}({\mathbf{x}})$ for some ${\mathbf{b}}^{\omega}\in V$. It is easy
to see that this defines another action of $\mathrm{Sp}_{2g}({\mathbb{F}}%
_{2})$ on the set $V$ given by $\omega\cdot{\mathbf{b}} := {\mathbf{b}%
}^{\omega}$. (Note that this action does not respect the structure of $V$ as a
vector space. In particular, it is different from the natural action given by
matrix multiplication.) It is easy to verify that the restriction of this
action to $\mathrm{Sp}_{2g-2}({\mathbb{F}}_{2})<\mathrm{Sp}_{2g}({\mathbb{F}%
}_{2})$ can be canonically identified with the analogous action of
$\mathrm{Sp}_{2g-2}({\mathbb{F}}_{2})$ on quadratic forms on ${\mathbb{F}}%
_{2}^{2g-2}$.

Now consider the two orbits ${\mathcal{Q}}^{+}$ and ${\mathcal{Q}}^{-}$ of
quadratic forms of type $+$ respectively type $-$ under the action of
$\mathrm{Sp}_{2g}({\mathbb{F}}_{2})$. We will show that ${\mathcal{Q}%
}^{\epsilon}$ ($\epsilon\in\{\pm\}$) splits up into four orbits under the
action of $\mathrm{Sp}_{2g-2}({\mathbb{F}}_{2})$, one for each of the $2^{2}$
possible values of $(b_{1},b_{2g})$, and that three of these orbits have the
type $\epsilon$ while the remaining orbit has the type $-\epsilon$.

Consider first the orbit ${\mathcal{Q}}^{+}$ of quadratic forms of type $+$,
that is, of Witt index $g$, and define
\[
{\mathbf{b}}_{0} = \left(
\begin{array}
[c]{c}%
0\\
0\\
0\\[-1ex]%
\vdots\\
0\\
0\\
0
\end{array}
\right)  \;, \quad{\mathbf{b}}_{1} = \left(
\begin{array}
[c]{c}%
1\\
0\\
0\\[-1ex]%
\vdots\\
0\\
0\\
0
\end{array}
\right)  \;, \quad{\mathbf{b}}_{2} = \left(
\begin{array}
[c]{c}%
0\\
0\\
0\\[-1ex]%
\vdots\\
0\\
0\\
1
\end{array}
\right)  \;, \quad{\mathbf{b}}_{3} = \left(
\begin{array}
[c]{c}%
1\\
1\\
0\\[-1ex]%
\vdots\\
0\\
1\\
1
\end{array}
\right)  \;.
\]

The quadratic forms $Q_{{\mathbf{b}}_{i}}({\mathbf{x}})$ for $i=0,1,2,3$ have
Witt index $g$, that is, are elements of ${\mathcal{Q}}^{+}$. Moreover,
${\mathcal{Q}}^{+} = \bigsqcup_{i=1}^{4} Q_{{\mathbf{b}}_{i}}^{\mathrm{Sp}%
_{2g-2}({\mathbb{F}}_{2})} \;. $ It is easy to see that the restriction of
$Q_{{\mathbf{b}}_{i}}$ to ${\mathbb{F}}_{2}^{2g-2}$ is of type $+$ (Witt index
$g-1$) for $i=0,1,2$ and of type $-$ (Witt index $g-2$) for $i=3$. Hence
$\overline{\phi}_{g}^{\,+}|_{\mathrm{Sp}_{2g-2}({\mathbb{F}}_{2})}%
\cong(\overline{\phi}_{g-1}^{\,+})^{3}\oplus\overline{\,\phi}_{g-1}^{-}$ as claimed.

The argument for the orbit ${\mathcal{Q}}^{-}$ of quadratic forms of type $-$
is analogous. \qed

\begin{proposition}
\label{11107P8.1}

\begin{itemize}
\item[\textbf{(1)}] Let $g\geq3$ and $n\geq1$. Then $\phi_{g,n}^{-}%
:{\mathcal{M}}_{g,n}\rightarrow{\mathfrak{S}}_{N_{g}^{-}}$ is equivalent to an
extension of the representation $(\phi_{g-1,n}^{-})^{3}\oplus\phi_{g-1,n}^{+}$
from ${\mathcal{M}}_{g-1,n}$ to ${\mathcal{M}}_{g,n}$, and $\phi_{g,n}%
^{+}:{\mathcal{M}}_{g,n}\rightarrow{\mathfrak{S}}_{N_{g}^{+}}$ is equivalent
to an extension of the representation $\phi_{g-1,n}^{-}\oplus(\phi_{g-1,n}%
^{+})^{3}$ from ${\mathcal{M}}_{g-1,n}$ to ${\mathcal{M}}_{g,n}$.

\item[\textbf{(2)}] Let $g\geq3$ and $n\geq0$. Let $b$ be a nonseparating
simple closed curve on $\Sigma_{g,n}$, and let $T_{b}$ be the Dehn twist
around $b$. Then the cycle structure of the image of $T_{b}$ under $\phi
_{g,n}^{-}$ is
\[
(1)^{2^{2g-2}}(2)^{2^{g-2}(2^{g-1}-1)}\,,
\]
and the cycle structure of the image of $T_{b}$ under $\phi_{g,n}^{+}$ is
\[
(1)^{2^{2g-2}}(2)^{2^{g-2}(2^{g-1}+1)}\,.
\]

\end{itemize}
\end{proposition}

\bigskip\noindent\textbf{Proof.} We first consider Part (1). For a simple
closed curve $c$ on $\Sigma_{g,1}$, we denote by $[c]$ the class of $c$ in
$H_{1}(\Sigma_{g,1},{\mathbb{Z}})$. We consider the curves $u_{1},\dots
,u_{g},v_{1},\dots,v_{g}$ indicated in Figure 7.1, and choose $[u_{g}%
],[u_{g-1}],\dots,[u_{1}],[v_{1}],[v_{2}],\dots,[v_{g}]$ as basis elements for
$H_{1}(\Sigma_{g,1},{\mathbb{Z}})$. With respect to this ordering, the
bilinear form $\langle-,\,-\rangle$ yielding the algebraic intersection number
is given by the matrix
\[
M=\left[
\begin{array}
[c]{ccccc}%
0 & 0 & \cdots & 0 & 1\\
0 & 0 & \cdots & 1 & 0\\
\vdots & \vdots &  & \vdots & \vdots\\
0 & 1 & \cdots & 0 & 0\\
1 & 0 & \cdots & 0 & 0
\end{array}
\right]  \,.
\]

\begin{figure}[tbh]
\bigskip\centerline{
\setlength{\unitlength}{0.5cm}
\begin{picture}(16.5,4)
\put(0,0){\includegraphics[width=8.25cm]{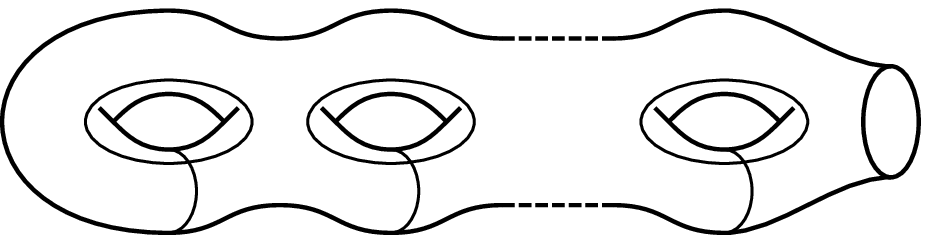}}
\put(2.7,0.5){\small $u_1$}
\put(6.7,0.5){\small $u_2$}
\put(12.7,0.5){\small $u_g$}
\put(3,3){\small $v_1$}
\put(7,3){\small $v_2$}
\put(13,3){\small $v_g$}
\end{picture}} \bigskip
\centerline{{\bf Figure 7.1.} Basis for $H_1(\Sigma_{g,1},\Z)$.}
\bigskip\end{figure}


\noindent The image of the Dehn twist $T_{b}$ about a simple closed curve $b$
under the epimorphism $\theta_{g,1}:{\mathcal{M}}_{g,1}\rightarrow
\mathrm{Sp}_{2g}({\mathbb{F}}_{2})$ is the map
\[
\lbrack a]\mapsto\lbrack a]+\langle a,b\rangle\lbrack b]\,.
\]
Since ${\mathcal{M}}_{g-1,1}<{\mathcal{M}}_{g,1}$ is generated by $T_{0}%
,T_{1},\dots,T_{2g-2}$ and since $a_{0},a_{1},\dots,a_{2g-2}$ intersect
neither $u_{g}$ nor $v_{g}$, it is clear from the above that $\theta
_{g,1}({\mathcal{M}}_{g-1,1})=\mathrm{Sp}_{2g-2}({\mathbb{F}}_{2}%
)<\mathrm{Sp}_{2g}({\mathbb{F}}_{2})$, and the restriction of $\theta_{g,1}$
to ${\mathcal{M}}_{g-1,1}$ coincides with $\theta_{g-1,1}:{\mathcal{M}%
}_{g-1,1}\rightarrow\mathrm{Sp}_{2g-2}({\mathbb{F}}_{2})<\mathrm{Sp}%
_{2g}({\mathbb{F}}_{2})$.

\medskip\noindent Let $n\geq1$. Gluing disks along all boundary components of
$\Sigma_{g,n}$ but one, we obtain an embedding $\Sigma_{g,n}\hookrightarrow
\Sigma_{g,1}$ that induces a surjective homomorphism $\mu_{g,n}:{\mathcal{M}%
}_{g,n}\twoheadrightarrow{\mathcal{M}}_{g,1}$ (see \cite{ParRol1}). More
precisely, the epimorphism $\mu_{g,n}$ is defined sending $T_{j}^{\prime}$ to
$T_{1}$ for all $j\in\{1,\dots,n\}$, and $T_{i}$ to $T_{i}$ for all
$i\in\{0,2,\dots,2g\}$. Observe that $\theta_{g,n}=\theta_{g,1}\circ\mu_{g,n}$
and the following diagram commutes
\[
\xymatrix{
\MM_{g-1,n} \ar[r] \ar[d]_{\mu_{g-1,n}} & \MM_{g,n} \ar[d]^{\mu_{g,n}}\\
\MM_{g-1,1} \ar[r] & \MM_{g,1}}
\]
where ${\mathcal{M}}_{g-1,n}\rightarrow{\mathcal{M}}_{g,n}$ and ${\mathcal{M}%
}_{g-1,1}\rightarrow{\mathcal{M}}_{g,1}$ are the natural embeddings described
in Section~\ref{11107S2}. By the above, it follows that $\theta_{g,n}%
({\mathcal{M}}_{g-1,n})=\mathrm{Sp}_{2g-2}({\mathbb{F}}_{2})<\mathrm{Sp}%
_{2g}({\mathbb{F}}_{2})$, and the restriction of $\theta_{g,n}$ to
${\mathcal{M}}_{g-1,n}$ coincides with $\theta_{g-1,n}:{\mathcal{M}}%
_{g-1,n}\rightarrow\mathrm{Sp}_{2g-2}({\mathbb{F}}_{2})<\mathrm{Sp}%
_{2g}({\mathbb{F}}_{2})$.

\medskip We have
\[
\overline{\phi}_{g}^{+}|_{\mathrm{Sp}_{2g-2}({\mathbb{F}}_{2})}\cong%
(\overline{\phi}_{g-1}^{+})^{3}\oplus\overline{\phi}_{g-1}^{-}\quad
\text{and}\quad\overline{\phi}_{g}^{-}|_{\mathrm{Sp}_{2g-2}({\mathbb{F}}_{2}%
)}\cong(\overline{\phi}_{g-1}^{-})^{3}\oplus\overline{\phi}_{g-1}^{+}%
\]
from Lemma~\ref{L:SymplecticPresentations}. By composing with $\theta_{g,n}$
we conclude that
\[
\phi_{g,n}^{+}|_{{\mathcal{M}}_{g-1,n}}\cong(\phi_{g-1,n}^{+})^{3}\oplus
\phi_{g-1,n}^{-}\quad\text{and}\quad\phi_{g,n}^{-}|_{{\mathcal{M}}_{g-1,n}%
}\cong(\phi_{g-1,n}^{-})^{3}\oplus\phi_{g-1,n}^{+}\,.
\]

\noindent Part (2) for $n \ge1$ follows by induction on $g$, using the cycle
structures $(2)^{6}\,(1)^{16}$ of $\phi_{3,n}^{-}(T_{b})$ and $(2)^{10}%
\,(1)^{16}$ of $\phi_{3,n}^{+}(T_{b})$ for $g=3$ (see Lemma \ref{11107L4.2}),
and the equivalence $\phi_{g,n}^{\epsilon}|_{{\mathcal{M}}_{g-1,n}} \cong%
(\phi_{g-1,n}^{\epsilon})^{3}\oplus\phi_{g-1,n}^{-\epsilon}$ proved above for
$g \ge4$.

\medskip Gluing a disk along the boundary component of $\Sigma_{g,1}$ we
obtain an embedding $\Sigma_{g,1}\hookrightarrow\Sigma_{g,0}$ that induces a
surjective homomorphism $\nu:{\mathcal{M}}_{g,1}\rightarrow{\mathcal{M}}%
_{g,0}$ (see \textsl{e.g. }\cite{ParRol1}). Observe that $\theta_{g,1}%
=\theta_{g,0}\circ\nu$. Thus, if $b$ is a nonseparating simple closed curve in
$\Sigma_{g,1}$, then $\theta_{g,1}(T_{b})=\theta_{g,0}(T_{b})$; therefore
$\phi_{g,1}^{\epsilon}(T_{b})=\phi_{g,0}^{\epsilon}(T_{b})$. This proves Part
(2) for $n=0$. \qed


\section{Surfaces with multiple boundary components\label{11107S9}}

We start the section with a geometrical interpretation of the decomposition
given in Proposition~\ref{11107P8.1}.

\bigskip\noindent Take $\epsilon\in\{ \pm\}$ and consider the representation
$\phi_{g,1}^{\epsilon}: {\mathcal{M}}_{g,1} \to{\mathfrak{S}}_{N_{g}%
^{\epsilon}}$. Let $w_{i}$ denote the image of $T_{i}$ under $\phi
_{g,1}^{\epsilon}$ for all $i \in\{0,1, \dots, 2g+1\}$. Set $H = \langle
w_{2g}, w_{2g+1} \rangle$. Since $w_{i}$ has order $2$, there is an
epimorphism from the group
\[
\widehat H =\langle x_{0},x_{1} \mid x_{0}^{2} = x_{1}^{2}=1,\ x_{0}x_{1}x_{0}
= x_{1} x_{0} x_{1} \rangle={\mathfrak{S}}_{3}
\]
to $H$ which sends $x_{i}$ to $w_{2g+i}$ for $i \in\{0,1\}$. In particular,
the order of $H$ divides $6$.

\medskip For $k\geq1$, let $\widehat{\Omega}_{k}$ denote the union of the
$k$-orbits $\Omega_{k,i}$ of $H$, $1\leq i\leq h_{k}$. Since the image
$\phi_{g,1}^{\epsilon}({\mathcal{M}}_{g-1,1})=\langle w_{0},\dots
,w_{2g-2}\rangle$ belongs to the centralizer of $H$, it acts on the set of
orbits $\{\Omega_{k,1},\dots,\Omega_{k,h_{k}}\}$ as well as on the whole set
$\widehat{\Omega}_{k}$. The first action induces a homomorphism $\nu
_{k}:{\mathcal{M}}_{g-1,1}\rightarrow{\mathfrak{S}}_{h_{k}}$ and the second a
homomorphism $\psi_{k}:{\mathcal{M}}_{g-1,1}\rightarrow{\mathfrak{S}}_{kh_{k}%
}$. The following lemma describes these representations.

\begin{lemma}
\label{11107L9.1} Let $g\geq3$. With the above notation, the following hold.

\begin{itemize}
\item[\textbf{(1)}] $h_{3}=N_{g-1}^{\epsilon}$, $h_{1}=N_{g-1}^{-\epsilon}$,
and $h_{k}=0$ for all $k \not \in \{1,3\}$.

\item[\textbf{(2)}] $\nu_{3} \cong\phi_{g-1,1}^{\epsilon}$, $\psi_{3}
\cong(\phi_{g-1,1}^{\epsilon})^{3}$, and $\nu_{1}=\psi_{1} \cong\phi
_{g-1,1}^{-\epsilon}$.

\item[\textbf{(3)}] $\{1, \dots, N_{g}^{\epsilon}\} = S(\phi_{g,1}^{\epsilon
}({\mathcal{M}}_{3,1}))$, where, for $G < {\mathfrak{S}}_{N_{g}^{\epsilon}}$,
$S(G)$ denotes the support of $G$, and ${\mathcal{M}}_{3,1}$ is the subgroup
of ${\mathcal{M}}_{g,1}$ generated by $T_{0},T_{1}, \dots, T_{6}$.
\end{itemize}
\end{lemma}

\medskip\noindent\textbf{Proof.} We prove by induction on $g$ that $\langle
w_{1},w_{2}\rangle$ has $N_{g-1}^{\epsilon}$ $3$-orbits, has no other
nontrivial orbit, and has $N_{g-1}^{-\epsilon}$ fixed letters. The case $g=3$
follows from Lemma \ref{11107L4.1}, and the inductive step can easily be
derived from the formula $\phi_{g,1}^{\epsilon}|_{{\mathcal{M}}_{g-1,1}}%
\cong(\phi_{g-1,1}^{\epsilon})^{3}\oplus\phi_{g-1,1}^{-\epsilon}$ of
Proposition ~\ref{11107P8.1}. Thanks to Lemma \ref{11107L2.2} this proves Part (1).

\medskip Let $\Omega_{3,i}$ be a $3$-orbit of $H$. We may write $\Omega
_{3,i}=\{a_{i},b_{i},c_{i}\}$ so that the restriction of $w_{2g}$ to
$\Omega_{3,i}$ is the transposition $(a_{i}\;b_{i})$ and the restriction of
$w_{2g+1}$ is $(b_{i}\;c_{i})$. Set $\widehat{\Omega}_{3}^{1}=\widehat{\Omega
}_{3}-S(w_{2g+1})=\{a_{i};1\leq i\leq N_{g-1}^{\epsilon}\}$, $\widehat{\Omega
}_{3}^{3}=\widehat{\Omega}_{3}-S(w_{2g})=\{c_{i};1\leq i\leq N_{g-1}%
^{\epsilon}\}$, and $\widehat{\Omega}_{3}^{2}=\widehat{\Omega}_{3}%
-(\widehat{\Omega}_{3}^{1}\cup\widehat{\Omega}_{3}^{3})=\{b_{i};1\leq i\leq
N_{g-1}^{\epsilon}\}$. Then $\widehat{\Omega}_{3}^{\ell}$ is invariant under
the action of ${\mathcal{M}}_{g-1,1}$ for $\ell=1,2,3$, and the action of
${\mathcal{M}}_{g-1,1}$ on $\widehat{\Omega}_{3}^{\ell}$ is equivalent to
$\nu_{3}$. So,
\[
\phi_{g,1}^{\epsilon}|_{{\mathcal{M}}_{g-1,1}}\cong\psi_{3}\oplus\psi_{1}%
\cong(\nu_{3})^{3}\oplus\nu_{1}\,.
\]
For the case $g=3$, one can easily check the equivalences $\nu_{3}\cong%
\phi_{2,1}^{\epsilon}$ and $\nu_{1}\cong\phi_{2,1}^{-\epsilon}$ using
Lemma~\ref{11107L3.1} and Lemma~\ref{11107L4.1}, so Part (2) holds in this
case and we can assume $g\geq4$. Recall that $h_{3}=N_{g-1}^{\epsilon}$ and
$h_{1}=N_{g-1}^{-\epsilon}$. By Lemma~\ref{11107L4.1} (if $g=4$), respectively
Proposition~\ref{11107P7.3} (if $g>4$), $\nu_{3}$ and $\nu_{1}$ are either
trivial or of the form $\nu_{3}=\phi_{g-1,1}^{\mu}\oplus\mathbf{1}_{q}$,
respectively $\nu_{1}=\phi_{g-1,1}^{\mu^{\prime}}\oplus\mathbf{1}_{q^{\prime}%
}$. Now, the formula of Proposition~\ref{11107P8.1} implies that $\phi
_{g,1}^{\epsilon}|_{{\mathcal{M}}_{g-1,1}}$ has no fixed letter, and thus both
$\nu_{3}$ and $\nu_{1}$ must be nontrivial and we must have $q=q^{\prime}=0$,
$\mu=\epsilon$, and $\mu^{\prime}=-\epsilon$. So Part (2) holds.

\medskip Part (3) follows by induction on $g$ using Lemma \ref{11107L4.1} for
the case $g=3$ and the equivalence $\phi_{g,1}^{\epsilon}|_{{\mathcal{M}%
}_{g-1,1}}\cong(\phi_{g-1,1}^{\epsilon})^{3}\oplus\phi_{g-1,1}^{-\epsilon}$
for the inductive step. \qed

\bigskip Now, we finish the proof of Theorem \ref{11107T0.2} with the following.

\begin{proposition}
\label{11107P9.2}Let $g\geq4$ and $n\geq0$.

\begin{itemize}
\item[\textbf{(1)}] ${\mathcal{O}}_{g,n}^{-}$ is the unique subgroup of
${\mathcal{M}}_{g,n}$ of index $N_{g}^{-}=2^{g-1}(2^{g}-1)$, up to conjugation.

\item[\textbf{(2)}] ${\mathcal{O}}_{g,n}^{+}$ is the unique subgroup of
${\mathcal{M}}_{g,n}$ of index $N_{g}^{+}=2^{g-1}(2^{g}+1)$, up to conjugation.

\item[\textbf{(3)}] All the other subgroups of ${\mathcal{M}}_{g,n}$ are of
index at least $5N_{g-1}^{-} > N_{g}^{+}$.
\end{itemize}
\end{proposition}

\bigskip\noindent\textbf{Proof.} Recall that we are under the inductive
hypothesis stated at the beginning of Part \ref{PartInduction}, that is,
Theorem \ref{11107T0.2} holds for a surface of genus $g-1$.

\medskip\noindent The case $n=1$ is proved in Proposition \ref{11107P7.3}, and
the case $n=0$ follows from Proposition \ref{11107P7.3}, Theorem
\ref{11107T0.1}, and the existence of the epimorphisms ${\mathcal{M}}%
_{g,1}\twoheadrightarrow{\mathcal{M}}_{g,0}$, ${\mathcal{M}}_{g,0}%
\twoheadrightarrow\mathrm{Sp}_{2g}({\mathbb{F}}_{2})$ described in Section
\ref{11107S0}. So, we may assume $n\geq2$.

\medskip\noindent Let $\varphi:{\mathcal{M}}_{g,n}\rightarrow{\mathfrak{S}%
}_{m}$ be a nontrivial transitive homomorphism with $m<5N_{g-1}^{-}$. As ever,
for $i\in\{0,2,\dots,2g+1\}$ and $j\in\{1,\dots,n\}$, we denote by $a_{i}$ and
$b_{j}$ the simple closed curves illustrated in Figure 2.2, we denote by
$T_{i}$ the Dehn twist about $a_{i}$, by $T_{j}^{\prime}$ the Dehn twist about
$b_{j}$, and we set $w_{i}=\varphi(T_{i})$ and $w_{j}^{\prime}=\varphi
(T_{j}^{\prime})$. For $j\in\{1,\dots,n\}$, we denote by $\Sigma^{(j)}$ a
tubular neighborhood of $b_{j}\cup a_{0}\cup(\cup_{i=2}^{2g}a_{i})$. Observe
that $\Sigma^{(j)}$ is a subsurface of $\Sigma_{g,n}$ of genus $g$ with a
unique boundary component, and the inclusion $\Sigma^{(j)}\hookrightarrow
\Sigma_{g,n}$ induces an injective homomorphism $\gamma_{j}:{\mathcal{M}%
}_{g,1}\rightarrow{\mathcal{M}}_{g,n}$ which sends $T_{1}$ to $T_{j}^{\prime}$
and $T_{i}$ to $T_{i}$ for all $i\in\{0,2,\dots,2g\}$. Set $\varphi
_{j}=\varphi\circ\gamma_{j}:{\mathcal{M}}_{g,1}\rightarrow{\mathfrak{S}}_{m}$.
By Proposition \ref{11107P7.3}, $\varphi_{j}$ is of the form $\varphi_{j}%
=\psi_{j}\oplus\mathbf{1}_{q_{j}}$ where $\psi_{j}$ is conjugate to an element
of $\{\phi_{g,1}^{+},\phi_{g,1}^{-}\}$ and $\mathbf{1}_{q_{j}}:{\mathcal{M}%
}_{g,1}\rightarrow{\mathfrak{S}}_{q_{j}}$ is the trivial representation.

\noindent Let $\Sigma^{\prime}$ be a tubular neighborhood of $\cup_{i=0}%
^{5}a_{2g-i}$. Then $\Sigma^{\prime}$ is a subsurface of $\Sigma_{g,n}$ of
genus $3$ with a unique boundary component, and is included in $\Sigma^{(j)}$
for all $j\in\{1,\dots,n\}$. Moreover, the inclusion $\Sigma^{\prime
}\hookrightarrow\Sigma_{g,n}$ induces an embedding $\gamma^{\prime
}:{\mathcal{M}}_{3,1}\rightarrow{\mathcal{M}}_{g,n}$. Set $\varphi^{\prime
}=\varphi\circ\gamma^{\prime}:{\mathcal{M}}_{3,1}\rightarrow{\mathfrak{S}}%
_{m}$. We have $S(\varphi_{j}({\mathcal{M}}_{g,1}))=S(\varphi^{\prime
}({\mathcal{M}}_{3,1}))$ by Lemma \ref{11107L9.1}\thinspace(3) for all $j$,
and $\bigcup_{j=1}^{n}S(\varphi_{j}({\mathcal{M}}_{g,1}))=S(\varphi
({\mathcal{M}}_{g,n}))$, thus
\[
S(\varphi({\mathcal{M}}_{g,n}))=S(\varphi_{j}({\mathcal{M}}_{g,1}%
))=S(\varphi^{\prime}({\mathcal{M}}_{3,1}))\,.
\]
Since $\varphi$ is transitive, it follows that there exists $\epsilon\in
\{\pm\}$ such that $m=N_{g}^{\epsilon}$ and $\varphi_{j}$ is conjugate to
$\phi_{g,1}^{\epsilon}$ for all $j\in\{1,\dots,n\}$.

\noindent Set $H=\langle w_{2g},w_{2g+1}\rangle$; note that $H \subseteq
\bigcap_{j}\gamma_{j}({\mathcal{M}}_{g,1})$. For $k\geq1$ we denote by
$\widehat{\Omega}_{k}$ the union of the $k$-orbits $\Omega_{k,i}$ of $H$,
$1\leq i\leq h_{k}$. On the other hand, we assume that ${\mathcal{M}}_{g-1,n}$
is the subgroup of ${\mathcal{M}}_{g,n}$ generated by $T_{1}^{\prime}%
,\dots,T_{n}^{\prime},T_{0},T_{2},\dots,T_{2g-2}$. Since the image
$\varphi({\mathcal{M}}_{g-1,n})$ lies in the centralizer of $H$, it acts on
the set of orbits $\{\Omega_{k,1},\dots,\Omega_{k,h_{k}}\}$ as well as on the
whole set $\widehat{\Omega}_{k}$. The first action induces a homomorphism
$\nu_{k}:{\mathcal{M}}_{g-1,n}\rightarrow{\mathfrak{S}}_{h_{k}}$, and the
second induces a homomorphism $\psi_{k}:{\mathcal{M}}_{g-1,n}\rightarrow
{\mathfrak{S}}_{kh_{k}}$. Applying Lemma \ref{11107L9.1}\thinspace(1) to
$\varphi_{j}({\mathcal{M}}_{g,1})$ for any $j$, we get that $h_{3}%
=N_{g-1}^{\epsilon}$, $h_{1}=N_{g-1}^{-\epsilon}$, and $h_{k}=0$ if
$k\not \in \{1,3\}$.

\noindent As in the proof of Lemma \ref{11107L9.1}, we set $\widehat{\Omega
}_{3}^{1}=\widehat{\Omega}_{3}-S(w_{2g+1})$, $\widehat{\Omega}_{3}%
^{3}=\widehat{\Omega}_{3}-S(w_{2g})$, and $\widehat{\Omega}_{3}^{2}%
=\widehat{\Omega}_{3}-(\widehat{\Omega}_{3}^{1}\cup\widehat{\Omega}_{3}^{3})$.
Then $\widehat{\Omega}_{3}^{\ell}$ is invariant under the action of
${\mathcal{M}}_{g-1,n}$ for $\ell=1,2,3$, and the action of ${\mathcal{M}%
}_{g-1,n}$ on $\widehat{\Omega}_{3}^{\ell}$ is equivalent to $\nu_{3}$.
Hence,
\[
\varphi|_{{\mathcal{M}}_{g-1,n}}\cong\psi_{3}\oplus\psi_{1}\cong(\nu_{3}%
)^{3}\oplus\nu_{1}\,.
\]
Since $h_{3},h_{1}\leq N_{g-1}^{+}$, by induction we have $\nu_{3}%
(T_{j}^{\prime})=\nu_{3}(T_{1}^{\prime})$ and $\nu_{1}(T_{j}^{\prime})=\nu
_{1}(T_{1}^{\prime})$, thus $\varphi(T_{j}^{\prime})=\varphi(T_{1}^{\prime})$
for all $j\in\{1,\dots,n\}$. Since $\varphi_{1}$ is conjugate to $\phi
_{g,1}^{\epsilon}$, we conclude that $\varphi$ is conjugate to $\phi
_{g,n}^{\epsilon}$. \qed

\bigskip\noindent\textbf{Proof of Theorem \ref{11107T0.3}.} This follows from
Proposition \ref{11107P8.1} and Theorem \ref{11107T0.2}, proved above. \qed




\bigskip\bigskip\noindent\textbf{A.\ Jon Berrick,}

\smallskip\noindent Department of Mathematics, National University of

Singapore, 10 Lower Kent Ridge Road, Singapore 119076, SINGAPORE.

\smallskip\noindent E-mail: \texttt{berrick@math.nus.edu.sg}

\bigskip\noindent\textbf{Volker Gebhardt,}

\smallskip\noindent School of Computing \& Mathematics, University of Western
Sydney, Locked Bag 1797, Penrith NSW 2751, AUSTRALIA.

\smallskip\noindent E-mail: \texttt{v.gebhardt@uws.edu.au}

\bigskip\noindent\textbf{Luis Paris,}

\smallskip\noindent Université de Bourgogne, Institut de Mathématiques de
Bourgogne, UMR 5584 du CNRS, B.P. 47870, 21078 Dijon cedex, FRANCE.

\smallskip\noindent E-mail: \texttt{lparis@u-bourgogne.fr}


\begin{thebibliography}{99}                                                                                               %


\bibitem {AraSou1}{\small \textbf{J.\ Aramayona, J.\ Souto.}
\textit{Homomorphisms between mapping class groups.} Preprint.
arXiv:1011.1855. }

\bibitem {Aschb1}{\small \textbf{M.\ Aschbacher.} \textit{On the maximal
subgroups of the finite classical groups.} Invent.\ Math.\ \textbf{76} (1984),
no.\ 3, 469--514. }

\bibitem {BaMiSe1}{\small \textbf{H.\ Bass, J.\ Milnor, J.-P.\ Serre.}
\textit{Solution of the congruence subgroup problem for $\mathrm{SL}_{n}$
($n\geq3$) and $\mathrm{Sp}_{2n}$ ($n\geq2$).} Inst.\ Hautes Études
Sci.\ Publ.\ Math.\ \textbf{33} (1967), 59--137. }

\bibitem {Birma1}{\small \textbf{J.\,S.\ Birman.} \textit{On Siegel's modular
group.} Math.\ Ann.\ \textbf{191} (1971), 59--68. }

\bibitem {BoCaPl1}{\small \textbf{W.\ Bosma, J.\ Cannon, C.\ Playoust.}
\textit{The \textsc{Magma} algebra system I: The user language.}
J.\ Symb.\ Comput.\ \textbf{24} (1997), 235--265,
http://magma.maths.usyd.edu.au/magma. }

\bibitem {Brids1}{\small \textbf{M.\,R.\ Bridson.} \textit{Semisimple actions
of mapping class groups on $\mathrm{CAT}(0)$ spaces.} Geometry of Riemann
surfaces, 1--14, London Math.\ Soc.\ Lecture Note Ser., 368, Cambridge
Univ.\ Press, Cambridge, 2010. }

\bibitem {Brids2}{\small \textbf{M.\,R.\ Bridson.} \textit{On the dimension of
CAT(0) spaces where mapping class groups act.} J. Reine Angew. Math., to
appear.\ arXiv:0908.0690. }

\bibitem {CCNPW1}{\small \textbf{J.\,H.\ Conway, R.\,T.\ Curtis,
S.\,P.\ Norton, R.\,A. Parker, R.\,A.\ Wilson.} \textit{Atlas of finite
groups. Maximal subgroups and ordinary characters for simple groups. With
computational assistance from J.\,G.\ Thackray.} Oxford University Press,
Eynsham, 1985. }

\bibitem {Coope1}{\small \textbf{B.\,N.\ Cooperstein.} \textit{Minimal degree
for a permutation representation of a classical group.} Israel
J.\ Math.\ \textbf{30} (1978), no.\ 3, 213--235. }

\bibitem {CoxMos1}{\small \textbf{H.\thinspace S.\thinspace M.\ Coxeter,
W.\thinspace O.\thinspace J.\ Moser.} \textit{Generators and relations for
discrete groups.} Springer-Verlag, Berlin-Göttingen-Heidelberg, 1957. }

\bibitem {DicksonTAMS1908}{\small \textbf{L.\,E.\ Dickson.}
\textit{Representations of the general symmetric group as linear groups in
finite and infinite fields}. Trans.\ Amer.\ Math.\ Soc.\ \textbf{9} (1908),
121--148. }

\bibitem {Dye1}{\small \textbf{R.\,H.\ Dye.} \textit{Interrelations of
symplectic and orthogonal groups in characteristic two.} J.\ Algebra
\textbf{59} (1979), no.\ 1, 202--221. }

\bibitem {Farb1}{\small \textbf{B.\ Farb.} \textit{Some problems on mapping
class groups and moduli space.} Problems on mapping class groups and related
topics, 11--55, Proc.\ Sympos.\ Pure Math., 74, Amer.\ Math.\ Soc.,
Providence, RI, 2006. }

\bibitem {FarMar1}{\small \textbf{B.\ Farb, D.\ Margalit.} \textit{A Primer on
Mapping Class Groups.} Princeton University Press, to appear. }

\bibitem {Funar1}{\small \textbf{L.\ Funar.} \textit{Two questions on mappping
class groups.} Proc. Amer. Math. Soc. \textbf{139} (2011), no. 1, 375--382. }

\bibitem {Gross1}{\small \textbf{E.\,K.\ Grossman.} \textit{On the residual
finiteness of certain mapping class groups.} J.\ London Math.\ Soc.\ (2)
\textbf{9} (1974/75), 160--164. }

\bibitem {Harer1}{\small \textbf{J.\,L.\ Harer.} \textit{Stability of the
homology of the mapping class groups of orientable surfaces.} Ann.\ of
Math.\ (2) \textbf{121} (1985), no.\ 2, 215--249. }

\bibitem {Hurwi1}{\small \textbf{A.\ Hurwitz.} \textit{Über algebraische
Gebilde mit eindeutigen Transformationen in sich.} Math.\ Ann.\ \textbf{41}
(1893), 403--442. }

\bibitem {Kerko1}{\small \textbf{S.\,P.\ Kerckhoff.} \textit{The Nielsen
realization problem.} Ann.\ of Math.\ (2) \textbf{117} (1983), no.\ 2,
235--265. }

\bibitem {Korkm1}{\small \textbf{M.\ Korkmaz.} \textit{Low-dimensional
homology groups of mapping class groups: a survey.} Turkish
J.\ Math.\ \textbf{26} (2002), no.\ 1, 101--114. }

\bibitem {LabPar1}{\small \textbf{C.\ Labruère, L.\ Paris.}
\textit{Presentations for the punctured mapping class groups in terms of Artin
groups.} Algebr.\ Geom.\ Topol.\ \textbf{1} (2001), 73--114. }

\bibitem {Liebe1}{\small \textbf{M.\,W.\ Liebeck.} \textit{On the orders of
maximal subgroups of the finite classical groups.} Proc.\ London
Math.\ Soc.\ (3) \textbf{50} (1985), no.\ 3, 426--446. }

\bibitem {Matsu1}{\small \textbf{M.\ Matsumoto.} \textit{A presentation of
mapping class groups in terms of Artin groups and geometric monodromy of
singularities.} Math.\ Ann.\ \textbf{316} (2000), no.\ 3, 401--418. }

\bibitem {Mumfo1}{\small \textbf{D.\ Mumford.} \textit{Abelian quotients of
the Teichmüller modular group.} J.\ Analyse Math.\ \textbf{18} (1967),
227--244. }

\bibitem {Paris1}{\small \textbf{L.\ Paris.} \textit{Small index subgroups of
the mapping class group.} J.\ Group Theory \textbf{13} (2010), no.\ 4,
613--618. }

\bibitem {ParRol1}{\small \textbf{L.\ Paris, D.\ Rolfsen.} \textit{Geometric
subgroups of mapping class groups.} J.\ Reine Angew.\ Math.\ \textbf{521}
(2000), 47--83. }

\bibitem {Powel1}{\small \textbf{J.\ Powell.} \textit{Two theorems on the
mapping class group of a surface.} Proc.\ Amer.\ Math.\ Soc.\ \textbf{68}
(1978), no.\ 3, 347--350. }

\bibitem {Putma1}{\small \textbf{A.\ Putman.} \textit{Cutting and pasting in
the Torelli group.} Geom.\ Topol.\ \textbf{11} (2007), 829--865. }

\bibitem {Sims1}{\small \textbf{C.\,C.\ Sims.} \textit{Computation with
finitely presented groups.} Encyclopedia of Mathematics and its Applications,
48. Cambridge University Press, Cambridge, 1994. }

\bibitem {Taylo1}{\small \textbf{D.\,E.\ Taylor.} \textit{The geometry of the
classical groups.} Sigma Series in Pure Mathematics, 9. Heldermann Verlag,
Berlin, 1992. }

\bibitem {Tits1}{\small \textbf{J.\ Tits.} \textit{Systèmes générateurs de
groupes de congruence.} C.\ R.\ Acad.\ Sci.\ Paris Sér.\ A-B \textbf{283}
(1976), no.\ 9, Ai, A693--A695. }

\bibitem {Wilso1}{\small \textbf{R.\,A.\ Wilson.} \textit{The finite simple
groups.} Graduate Texts in Mathematics, 251. Springer-Verlag London, Ltd.,
London, 2009. }

\bibitem {Zimme1}{\small \textbf{B.\,P.\ Zimmermann.} \textit{A note on
minimal finite quotients of mapping class groups.} Preprint. arXiv:0803.3144.
}
\end{thebibliography}
\end{document}